\documentclass[12pt]{article}

\usepackage{amsmath}
\usepackage{amsthm}
\usepackage[mathscr]{eucal}
\usepackage{amssymb,latexsym}
\usepackage{psfrag}
\usepackage{graphicx}

\usepackage{version} 

\theoremstyle{plain}
\newtheorem{thm}{Theorem}[section]
\newtheorem{lem}[thm]{Lemma}%[section]
\newtheorem{Def}{Definition}[section]
\newtheorem{prop}[thm]{Proposition}%[section]
\newtheorem{cor}[thm]{Corollary}%[section]
\newtheorem*{claim}{Claim} 

\newcommand\calA{{\mathcal{A}}}
\newcommand\calB{{\mathcal{B}}}

\newcommand\calF{{\mathcal{F}}}

\newcommand\calV{{\mathcal{V}}}
\newcommand{\tr}{\mbox{tr}}
\newcommand\calH{{\mathcal{H}}}

\newcommand\bbR{{\mathbb R}}

\newcommand\bbH{{\mathbb{H}}}

\renewcommand\l{\lambda}

\renewcommand\S{\Sigma}

\renewcommand\d{\partial}

\newcommand\f{\phi}
\renewcommand\L{\triangle}
\newcommand\D{\nabla}
\newcommand\e{\epsilon}
\renewcommand\b{\beta}

\renewcommand\div{{\rm div}}

\newcommand\la{\langle}
\newcommand\ra{\rangle}

\newcommand\ric{{\rm Ric}}
\renewcommand\l{\lambda}
\newcommand\g{\gamma}

\renewcommand\a{\alpha}

\newcommand\<{\la}
\renewcommand\>{\ra}

\newcommand{\half}{\frac{1}{2}} 
\newcommand{\muup}{\bar\mu}
\newcommand{\mulow}{\underline{\mu}}
\renewcommand\Re{{\mathbb R}}
 
\newcommand{\Ric}{\text{\rm Ric}}
\newcommand{\hRic}{\widehat{\Ric}}
\newcommand{\bu}{\bar u} 
\newcommand{\hS}{\hat S} 
\newcommand{\eps}{\epsilon}

\newcommand\beq{\begin{equation}}
\newcommand\eeq{\end{equation}}
\newcommand\ben{\begin{enumerate}}
\newcommand\een{\end{enumerate}}
\newcommand\bit{\begin{itemize}}
\newcommand\eit{\end{itemize}}

\newcommand{\tM}{\tilde M}
\newcommand{\tme}{\tilde g} 
\newcommand{\tS}{\tilde S} 
\newcommand{\hme}{\hat g}
\newcommand{\htme}{\hat{\tme}}
\newcommand{\that}{{\hat t}}

\newcommand{\hgamma}{{\hat \gamma}}

\newcommand{\arcsinh}{\text{arcsinh}}

\newcommand{\tRic}{\widetilde{\Ric}}
\newcommand{\tnabla}{\widetilde{\nabla}}

\newcommand{\tGamma}{\widetilde{\Gamma}}

\newcommand{\hpara}{h}
\newcommand{\nablapara}{{\nabla\!\!\!\!/}}
\newcommand{\Ricpara}{\text{ric}}
\newcommand{\Spara}{\text{s}}

\newcounter{mnotecount}[section]

\includeversion{ver:rigproof} 
\title{Rigidity and Positivity of Mass for Asymptotically Hyperbolic Manifolds}

\author{
L.  Andersson$^{a,b}$, M.Cai$^{a}$, G. J. Galloway$^a$
\\ \\
${}^a$Department of Mathematics, University of Miami \\ Coral Gables, FL
33124
\\ \\
${}^b$Max-Planck-Institut f\"ur Gravitationsphysik\\
Am M\"uhlenberg 1, 14476 Golm,
Germany
}

\begin{document}

\date{}
%\date{\today \ {\em File:\jobname{.tex}}}
\maketitle

\begin{abstract}
The Witten spinorial argument has been adapted in several works
over the years to prove positivity of mass in the asymptotically AdS
and asymptotically hyperbolic settings in arbitrary dimensions.  
In this paper we prove a scalar curvature rigidity result and a 
positive mass theorem for asymptotically
hyperbolic manifolds that do not  require a spin assumption. The positive mass
theorem is reduced to the rigidity case by a deformation 
construction near the conformal boundary.  The proof of the rigidity 
result is based on a study of minimizers of the BPS brane action.
\end{abstract}

\section{Introduction}
Developments in string theory during the past decade, in particular the
emergence of the AdS/CFT correspondence, 
have increased
interest in the mathematical and physical properties of asymptotically
hyperbolic Riemannian manifolds. Such manifolds arise naturally as spacelike
hypersurfaces in asymptotically anti-de Sitter spacetimes. 

Asymptotically hyperbolic  manifolds have a rich geometry at
infinity, as exhibited by e.g., renormalized volume and $Q$-curvature. 
The mass of an asymptotically hyperbolic manifold may, 
under suitable asymptotic conditions, be defined as the
integral of a function defined at infinity, the so-called mass aspect
function. This feature is related to the fact that, in contrast to the
asymptotically Euclidean case, harmonic functions on an
AH manifold are not in general constant, but have
nontrivial boundary values at infinity. 

In this paper we shall prove a scalar curvature rigidity result and a 
positive mass theorem for asymptotically
hyperbolic manifolds.
The results do not  require a spin assumption. The positive mass
theorem is reduced to the rigidity case by a novel deformation 
construction near
the conformal boundary.  The proof of the rigidity result is based on 
a study of minimizers of the BPS brane action.

Let
$(M^{n+1},g)$ denote an $(n+1)$-dimensional Riemannian manifold, and let
$\mathbb H^{n+1}$ denote $(n+1)$-dimensional hyperbolic space of curvature $K
= -1$.  
As a prelude to proving positivity of mass
in the asymptotically hyperbolic setting 
(see the discussion below),
we  first establish the following rigidity result.
\begin{thm}\label{rigid}
Suppose $(M^{n+1},g)$, $2 \le n \le 6$, has scalar curvature $S[g]$ 
satisfying, $S[g]\ge -n(n+1)$,
and is isometric to $\mathbb H^{n+1}$  outside a compact set.
Then $(M^{n+1},g)$ is globally isometric to $\mathbb H^{n+1}$.
\end{thm}

In the case that $(M^{n+1},g)$ is a spin manifold, this theorem 
follows from
a result of Min-oo \cite{minoo} (see also \cite{AD,delay}), as well as from  the
rigidity part of the more recently proved positive mass theorem for
asymptotically hyperbolic manifolds  \cite{wang,  ChH}.  The main point of
Theorem \ref{rigid} is that it does not require 
a spin assumption.  We note,
for comparison, that there have been some other recently obtained rigidity
results for hyperbolic space \cite{qing, bonini:miao:qing:ricci,shi,cai:qing} 
that do not require a spin assumption, but
these impose conditions on the Ricci curvature.

The proof of Theorem \ref{rigid} is based on the general minimal surface
methodology of Schoen and Yau \cite{SY}, adapted to a negative lower bound on
the scalar curvature.    This means, in our approach, that  minimal surfaces
are replaced by non-zero constant mean curvature surfaces, and the area
functional is replaced by the so-called BPS brane action, as utilized by
Witten and Yau \cite{WY} in their work on the AdS/CFT correspondence.  
From
the regularity results of geometric measure theory, we require $M$ to have
dimension $\le 7$ in order to avoid the occurrence of singularities in
co-dimension one minimizers of the brane action
\footnote{\samepage However,
the work of Christ and Lohkamp \cite{christ:loh, loh2} offers the possibility of eliminating this dimension
restriction.}.  In Section 2 we prove a
local warped product splitting result, where the splitting takes place about
a certain minimizer of the brane action.   This splitting result, which
extends to the case of negative lower bound on the scalar curvature previous
results of Cai and Galloway \cite{CG, cai}, is then used to prove
Theorem~\ref{rigid}.

Our original motivation for proving Theorem \ref{rigid} was to obtain a proof
of positivity of mass for asymptotically hyperbolic manifolds that does not
require a spin assumption.   In \cite{gibbonsetal}, Gibbons et al. adapted
Witten's spinorial argument to prove positivity of mass in the $3+1$ asymptotically
AdS setting.  More recently, Wang \cite{wang}, and, under  weaker asymptotic conditions,
Chru\'sciel and Herzlich~\cite{ChH} have provided precise definitions of the mass in the asymptotically hyperbolic setting and have given spinor based proofs of positivity of mass in dimensions $\ge 3$.  These latter positive mass results may be
paraphrased  as follows:

\begin{thm}\label{pmass}
Suppose $(M^{n+1},g)$, $n \ge 2$, is an asymptotically hyperbolic spin manifold with
scalar curvature $S\ge -n(n+1)$.  Then $M$ has mass $m \ge 0$, and $=0$ iff M
is isometric to standard hyperbolic space $\mathbb H^{n+1}$.
\end{thm}

Physically, $M$ corresponds to a maximal (mean curvature zero)  spacelike
hypersurface in spacetime satisfying the Einstein equations with cosmological
constant $\Lambda = -n(n+1)/2$.  For then the Gauss equation and weak energy
condition imply $S \ge -n(n+1)$.

Here we present the following version of Theorem \ref{pmass}, which does not require $M$ to be spin.

\begin{thm} \label{thm:posmass} 
Let $(M^{n+1},g)$, $2 \leq n \leq 6$, be an asymptotically hyperbolic
manifold with scalar curvature $S[g] \ge -n(n+1)$. Assume that the mass
aspect function does not change sign, i.e. that it is either negative,
zero, or positive. Then, either the mass of $(M,g)$ is positive, or $(M,g)$
is isometric to hyperbolic space. 
\end{thm}
As noted above, the mass aspect function is a scalar function  whose integral over conformal infinity
determines the mass; see Section 3 for precise definitions.

 Our approach to proving Theorem \ref{thm:posmass} is inspired by Lohkamp's variation \cite{loh} of the 
Schoen-Yau~\cite{SY1} proof of the classical positive mass theorem for asymptotically
flat manifolds.  Our proof  makes use of Theorem \ref{rigid},  together with a
deformation result, which shows roughly that if an asymptotically hyperbolic
manifold with scalar curvature satisfying, $S\ge -n(n+1)$, has negative mass aspect
then the metric can be deformed near infinity to the hyperbolic metric, while
maintaining the scalar curvature inequality.  This deformation result (Theorem 
\ref{defm}), along with an analysis of the case in which the mass aspect
vanishes identically (Theorem \ref{thm:rigid2}), and their 
application to the proof of Theorem \ref{thm:posmass} are presented in Section 3. 

\section{The rigidity result}

The aim of this section is to give a proof of Theorem \ref{rigid}.

\subsection{The brane action}

Let $(M^{n+1}, g)$ be an $(n+1)$-dimensional oriented Riemannian manifold
with volume form $\Omega$.  Assume there is a globally defined
form $\Lambda$ such that $\Omega = d \Lambda$.

Let $\S^n$ be a compact orientable hypersurface in $M$.   Then $\S$ is
$2$-sided in $M$; designate one side as the ``outside" and the other as the
``inside".  Let $\nu$ be the {\it outward} pointing unit normal along $\S$, and
let $\S$ have the orientation induced by $\nu$ (i.e., determined by
the induced volume form $\omega = i_{\nu}\Omega$).   Then, for any such
$\S$, we define the brane action $\calB$ by,
\beq
\calB(\S) = \calA(\S) -n \calV(\S) \, ,
\eeq
where $\calA(\S) =$ the area of $\S$, and $\calV(\S) =  \int_{\S} \Lambda$.
If $\S$ bounds to the inside then,  by Stokes theorem, $\calV(\S) =$ the volume of the region enclosed by $\S$.   Although $\Lambda$ is not uniquely determined,
Stokes theorem shows that, within a given homology class, $\calB$ is uniquely
determined up to an additive constant.

We wish to consider the formulas for the first and second variation of the
brane action.  First, to fix notations, let $A$ denote the second fundamental
form of $\S$; by our conventions, for each pair of tangent vectors $X, Y \in T_p\S$,
\beq
A(X,Y) = \<\D_X\nu,Y\>  \, ,
\eeq
where $\D$ is the Levi-Civita connection of $(M,g)$ and $h = \<\,,\,\>$ is the
induced metric on $\S$. Then $H = {\rm tr}\, A$ is the mean curvature of $\S$.

Let $t \to \S_t$, $-\e <t < \e$, be a normal variation of $\S= \S_0$, with variation vector
field $V = \left . \frac{\d}{\d t}\right |_{t=0} = \phi \nu$,  $\phi \in C^{\infty}(\S)$.   Abusing notation slightly, set
$\calB(t) = \calB(\S_t)$.  Then for first variation we have,
\beq
\calB'(0) = \int_{\S}  (H -n)\f \, dA
\eeq
Thus $\S$ is a stationary point for the brane action if and only if it
has constant mean curvature $H = n$.

Assuming $\S$ has mean curvature
$H =n$, the second variation formula is given by
\beq
\calB''(0) = \int_{\S} \phi L(\phi)\, dA  \,,
\eeq
where,
\beq\label{op1}
L(\phi) = - \L \phi + \frac12(S_{\S} - S - |A|^2 - H^2)\,\phi   \,,
\eeq
and where $S_{\S}$ is the scalar curvature of $\S$ and $S$ is the scalar
curvature of $M$.  Here  $L$ is the {\it stability operator} associated with the brane action, and is closely related to the stability operator of minimal surface theory.
Using the fact that $H =n$,  $L$ can be re-expressed as,
\beq\label{op2}
L(\phi) = - \L \phi + \frac12(S_{\S} - S_n - |A_0|^2)\,\phi   \,,
\eeq
where $S_n = S + n(n+1)$ and $A_0$ is the trace free part of $A$, $A_0 = A -h$.
We note that, in our applications, $S_n$ will be nonnegative.

A stationary point $\S$ for the brane action is said to be $\calB$-stable provided
for all normal variations $t \to \S_t$ of $\S$, $\calB''(0) \ge 0$.
For operators of the form (\ref{op2}), the following proposition is well-known.
\begin{prop}\label{stable}
The following conditions are equivalent.
\ben
\item $\S$ is $\calB$-stable.
\item $\lambda_1 \ge 0$, where $\lambda_1$ is the principal eigenvalue of $L$.
\item There exists $\f \in C^{\infty}(\S)$, $\f > 0$, such that
$L(\phi) \ge 0$.
\een
\end{prop}
In particular, if $\l_1 \ge 0$, $\f$ in part 3 can be chosen to be an eigenfunction.

 \subsection{Warped product splitting}

 In this section we prove the local warped product splitting result alluded to in the
 introduction.  As a precursor, we prove the following infinitesimal rigidity result.

\begin{prop}\label{infrigid}
Let $(M^{n+1},g)$ be an oriented Riemannian manifold
with scalar curvature $S$ satisfying,
\beq\label{scalar}
S \ge -n(n+1) \, .
\eeq
Let $\S^n$ be a compact orientable
$\calB$-stable hypersurface in $M$ which does not admit a metric of positive
scalar curvature.   Then the following must  hold.
\ben
\item[(i)] $\S$ is umbilic, in fact $A =h$, where $h$ is the induced metric on $\S$.
\item[(ii)] $\S$ is Ricci flat and $S = -n(n+1)$ along $\S$.
\een
\end{prop}

\proof
By Proposition \ref{stable}, there exists $\f \in \S$, $\f > 0$, such that $L(\f) \ge 0.$
The scalar curvature $\tilde S$ of $\S$
in the conformally rescaled metric $\tilde h = \f^{\frac2{n-2}}h$ is then given by,
\begin{align}\label{rescale}
\tilde S & =  \phi^{-\frac{n}{n-2}}(-2\L\f + S_{\S} \f +\frac{n-1}{n-2}\frac{|\D\f|^2}{\f}) \nonumber \\\
& =  \phi^{-\frac{2}{n-2}}(2\f^{-1}L(\f) + S_n + |A_0|^2 +\frac{n-1}{n-2}\frac{|\D\f|^2}{\f^2})
\end{align}
where, for the second equation, we have used (\ref{op2}) with $f = \f$.  Since all terms
in parentheses above are nonnegative, (\ref{rescale}) implies that $\tilde S \ge 0$.
If $\tilde S > 0$ at some point, then by well known results
\cite{KW} one can conformally change $\tilde h$  to a metric
of strictly positive scalar curvature, contrary to assumption.  Thus $\tilde S$ vanishes
identically, which implies   $L(\f) =0$, $S_n = 0$, $A_0 = 0$ and $\f$ is constant.  Equation (\ref{op2}), with $f = \f$  then  implies that $S \equiv  0$.  By a result of Bourguinon (see \cite{KW}),
it follows that $\S$ carries a metric of positive scalar curvature unless  it is
Ricci flat.  Thus conditions (i) and (ii) are satisfied. ~\qed

Proposition \ref{rigid} will be used in the proof of the following local
warped product splitting result.

\begin{thm}\label{warp}
Let ($M^{n+1},g)$ be an oriented Riemannian manifold
with scalar curvature $S \ge -n(n+1)$.  Let $\S$ be a compact orientable
hypersurface in $M$ which does not admit a metric of positive scalar curvature.
If $\S$ locally minimizes the brane action $\calB$ then there is a neighborhood
$U$ of $\S$ such that $(U,g|_U)$ is isometric to the warped product
$((-\e,\e)\times \S, dt^2 +e^{2t} h)$, where $h$, the induced metric on $\S$, is
Ricci flat.
\end{thm}

By  ``locally minimizes" we mean, for example, that $\S$ has brane action less
than or equal to that of all graphs over $\S$ with respect to Gaussian normal
coordinates.  A related result has been obtained by Yau \cite{yau} 
in dimension three.

\smallskip
\proof
Let $\calH(u)$ denote the mean curvature of the hypersurface $\S_u : x \to exp_x u(x)\nu$, $u\in C^{\infty}(\S)$, $u$ sufficiently small.  $\calH$ has linearization $\calH'(0) = L$,
where $L$ is the $\calB$-stability operator~(\ref{op2}). But by Proposition \ref{infrigid},
$L$ reduces to $-\triangle$, and hence $\calH'(0) = - \triangle$.  We introduce the
operator,
\beq
\calH^* : C^{\infty}(\S) \times \bbR \to  C^{\infty}(\S) \times \bbR  \,, \quad
\calH^*(u,k) = \left(\calH(u) -k , \int_{\S}u\right)  \,,
\eeq
which one easily checks has invertible linearization about $(0,0)$, since the
kernel of $\calH'(0)$ contains only the constants.  By the inverse function theorem,
for each $\tau$ sufficiently small there exists $u = u_{\tau}$ and $k = k_{\tau}$
such that $\calH(u_{\tau}) = k_{\tau}$ and $\int_{\S} u_{\tau} dA = \tau$.   Since
$u'(0) \in {\rm ker}\, \calH'(0)$, the latter equation implies that $u'(0) =$ const $> 0$.
Thus for $\tau$ suffiiciently small, the hypersurfaces $\S_{u_{\tau}}$ form a foliation
of a neighborhood $U$ of $\S$ by constant mean curvature hypersurfaces.

Using coordinates on $\S$ and the normal field to the $\S_{u_{\tau}}$'s to transport these
coordinates to each $\S_{u_{\tau}}$, we have, up to isometry,
\beq\label{coords}
U = (-\e, \e) \times \S \qquad g|_U = \phi^2 dt^2 + h_t \,,
\eeq
where $h_t = h_{ij}(t,x)dx^idx^j$, $\phi = \phi(t,x)$ and $\S_t = \{t\} \times \S$
has constant mean curvature.   Since $\S$ locally minimizes the brane action,
we have, $\calB(0) \le \calB(t)$ for all $t \in (-\e,\e)$, for $\e$ sufficiently small.

Let $H(t)$ denote the
mean curvature of $\S_t$.    $H = H(t)$ obeys the evolution
equation,
\beq\label{evolve}
\frac{dH}{dt} = L(\f) \,,
\eeq
where for each $t$, $L$ is the operator on $\S_t$ given in Equation (\ref{op1}).
Since $\S$ locally minimizes the brane action, we have
$H(0) = n$.  We show $H \le n$ for $t \in [0,\e)$.  If this is not the case, there exists $t_0 \in (0,\e)$
such that $H(t_0) > n$.  Moreover, $t_0$ can be chosen so that $H'(t_0) > 0$.
Let $\tilde S$ be the scalar curvature of $\S_{t_0}$ in the conformally related
metric $\tilde h = \f^{\frac2{n-2}}h_{t_0}$.  Arguing similarly as in the derivation of
(\ref{rescale}), Equations (\ref{op1}) and (\ref{evolve}) imply,
\beq\label{rescale2}
\tilde S =  \phi^{-\frac{2}{n-2}}(2\f^{-1}H'(t_0)  + S + |A|^2 + H^2 +\frac{n-1}{n-2}\frac{|\D\f|^2}{\f^2})  \,,
\eeq
where all terms are evaluated on $\S_{t_0}$.  The Schwartz inequality gives, $|A|^2 \ge H^2/n > n$.  This, together with the assumed scalar curvature inequality (\ref{scalar}),
implies that $ S + |A|^2 + H^2 > 0$.  We conclude from (\ref{rescale2}) that
$\S_{t_0}$ carries a metric of positive scalar curvature, contrary to assumption.

Thus, $H \le n$ on $[0,\e)$, as claimed.  Now, by the formula for the first
variation of the brane action, it follows that
\beq\label{firstvar}
\calB'(t)  =  \int_{\S_t}  (H -n) \f\, dA \, \le 0 \,, \quad \mbox{for all  } t \in [0,\e) \,.
\eeq
But since $\calB$ achieves a minimum at $t=0$, it must be that
$\calB'(t)= 0$ for $t \in [0,\e)$. Hence,  the integral in (\ref{firstvar})
vanishes, which implies that $H = n$ on  $[0,\e)$.
A similar argument shows that $H = n$ on $(-\e, 0]$,
as well.  Equation (\ref{evolve}) then implies that $L(\f) = 0$ on each $\S_t$.
Hence, by Proposition \ref{stable}, each $\S_t$ is $\calB$-stable.  From
Proposition~\ref{infrigid}, we have that $A_t = h_t$, where $A_t$ is the second fundamental
form of $\S_t$, and that $\f$ only depends on $t$.  By a simple change of $t$-coordinate
in (\ref{coords}), we may assume without loss of generality that $\f = 1$.
Then the condition $A_t = h_t$ becomes, in the coordinates (\ref{coords}),
$\frac{\d h_{ij}}{\d t} = 2 h_{ij}$.  Upon integration this gives, $h_{ij}(t,x) = e^{2t}
h_{ij}(0,x)$, which completes the proof of the theorem. \qed

\subsection{Proof of the rigidity result}

In order to prove Theorem \ref{rigid} it is convenient to work with an
explicit representation of hyperbolic space $\bbH^{n+1}$.  We start with
the half-space model $(H^{n+1},g_H)$, where, $H^{n+1} = \{(y, x^1,\cdots, x^n) : y > 0\}$, and
\beq
g_H = \frac1{y^2} \left (dy^2 +  (dx^1)^2 + \cdots + (dx^n)^2 \right)  \,,
\eeq
and make the change of variable $y= e^{-t}$, to obtain
$\bbH^{n+1}= (\bbR^{n+1}, g_0)$,
where,
\beq
g_0 = dt^2 + e^{2t} \left( (dx^1)^2 + \cdots + (dx^n)^2 \right)  \,.
\eeq

As in the statement of Theorem \ref{rigid}, let $(M^{n+1},g)$ be a Riemannian manifold with scalar curvature $S[g]$ satisfying
$S[g] \ge -n(n+1)$.  We assume that there are compact sets $K \subset M$, $K_0
\subset \bbR^{n+1}$ such that  $M - K$ is diffeomorphic to $\bbR^{n+1} - K_0$, and,
with respect to Cartesian coordinates $(t,x_1,...,x_n)$ on the complement of $K$, $g = g_0$.   We want to show that $(M^{n+1},g)$ is globally
isometric to $\bbH^{n+1}$.  Since $M$ is simply connected near infinity, it is in fact
sufficient to show that $M$ has everywhere constant curvature $K_M = -1$.
Our approach to proving this is to partially compactify $(M,g)$ and then minimize the
brane action in a suitable homology class.

To partially compactify, we use the fact that the translations $x^i \to x^i + x^i_0$
are isometries on $(\bbR^{n+1}, g_0)$.  Choosing $a>0$ sufficiently large, we
can enclose the compact set $K$ in an infinitely long rectangular box, with sides
determined by the ``planes",
$x^i = \pm a, i = 1, \cdots , n, -\infty < t < \infty$.  We can then identify points on
opposite sides, $x^i = -a$, $x^i =a$, $i = 1, ...,n$, of the box in the obvious manner to obtain an identification space
which we denote by $(\hat M, \hat g)$.  Note that outside the compact set $K$,
\beq\label{cusp}
\hat M = \bbR \times T^n \,,  \quad \hat g = dt^2 + e^{2t} h \,,
\eeq
where $h$ is a flat metric on the torus $T^n$.  Thus, $(\hat M, \hat g)$ is just
a standard hyperbolic cusp outside the compact set $K$, 
with scalar curvature satisfying $S[\hat g] \ge -n(n+1)$
globally.

Choose $b > 0$ large so that $K$ is contained in
the region of $\hat M$ bounded between the toroidal slices $t = \pm b$,
and fix
a $t$-slice $\S_0 = \{t_0\} \times T^n$, $t_0 > b$.  $\S_0$ separates $\hat M$
into an ``inside" and an ``outside", the inside being the component of
$\hat M - \S_0$ containing the cusp end $t = -\infty$.  We consider the brane
action of hypersurfaces $\S$ homologous to $\S_0$,
\beq
\calB(\S) = A(\S) -n V(\S) \, .
\eeq
We note, as is needed to define $\calB(\S)$ unambiguously, that since $\S$ is homologous
to $\S_0$ it, too, has a distinguished ``inside" and ``outside", determined
by the
fact that both $\S_0$ and $\S$ are homologous to a $t$-slice far out on the
cusp end.

We now want to minimize the brane action $\calB$ in the homology class
$[\S_0]$.  The basic approach is to consider a minimizing sequence
$\S_1, \S_2, ...$ and use the compactness results of geometric measure
theory to extract a regular limit surface.   The potential difficulty with this
approach is that, in principle, the surfaces $\S_1, \S_2, ...$, or portions
of them, may drift out to infinity along either end of $\hat M$.  But, in fact, that can
be avoided in the present situation, owing to the existence of natural barrier
surfaces, namely the $t$-slices themselves.

Fix $t$-slices $\S' = \{t_1\} \times \S$, $t_1 > t_0$, and $\S'' = \{t_2\} \times \S$,
$t_2 < - b$.  We show that any minimizing sequence can be replaced by a minimizing
sequence contained in the region between $\S'$ and $\S''$.   To this end,
consider  a hypersurface $\S$ homologous to $\S_0$ that extends beyond $\S''$
into the region $t < -t_2$.  Without loss of generality we may assume $\S$ meets
$\S''$ transversely.  Let $D$ be the part of $\S$ meeting $\{t \le -t_2\}$, and
let $U$ be the domain bounded by $\S''$ and $D$.  Then $\d U$ consists
of $D$ and a part $D''$ of $\S''$.  Let $\hat \S$ be the hypersurface homologous
to $\S_0$ obtained from $\S$ by replacing $D$ with $D''$ (see Figure 2).
\begin{figure}[Ht]
\psfrag{XSS1}{$\S'$}
\psfrag{XSS0}{$\S_0$}
\psfrag{XSS2}{$\S''$}
\psfrag{XSSH}{$\hat \S$}
\psfrag{XSS}{$\S$}
\psfrag{XU}{$U$}
\psfrag{XK}{$K$}
\psfrag{XT}{$t$}
\centering 
\includegraphics[width=4.6in]{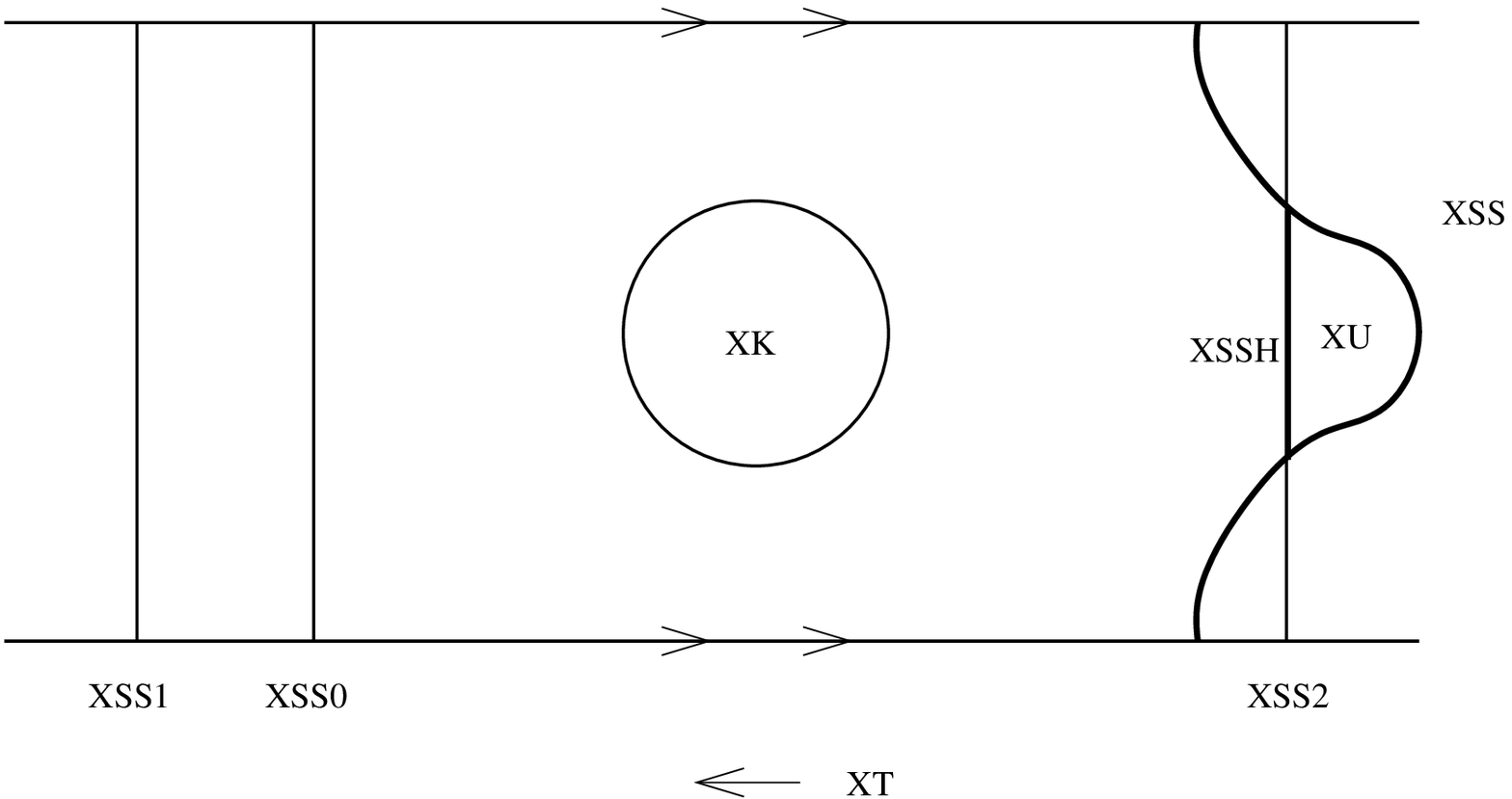} 
\caption{Replacing $\S$ by $\hat \S$}
\label{fig1}
\end{figure}

Since $U$ is  contained in a region where the metric ($\ref{cusp}$) applies,
and since in this region $\div(\d_t) = n$,  we apply the divergence
theorem to obtain,
\begin{align*}
n \calV(\hat \S) - n\calV(\S) & =  n \,{\rm vol}(U) = \int_U \div(\d_t) dV  \nonumber \\
& =  \int_D \<\d_t, n\> dA - \int_{D''} \<\d_t,\d_t \> dA  \nonumber \\
& \le  A(D) - A(D'')  =  A(\S) - A(\hat \S)  \,.
\end{align*}
Rearranging this inequality gives the  desired, $\calB(\hat\S) \le \calB(\S)$.
By a similar argument the same conclusion holds if $\S$ extends beyond
$\S'$ into the region $\{ t > t_1\}$.

Thus, we can choose a minimizing sequence $\S_i$, for the brane action
within the homology class $[\S_0]$ that is confined to the compact region
between $\S'$ and $\S''$.   Since $\calV(\S_i) \le \calV(\S')$, we are
ensured that $\lim_{i \to \infty} \calB(\S_i)  > - \infty$.  Then the compactness and
regularity results of geometric measure theory (see e.g., \cite{federer}, \cite{SY}
and references therein) guarantee the existence
of a regular embedded hypersurface $S$ homologous to $\S_0$ that
achieves a minimum of the brane action on $[\S_0]$.   In general
$S$ will be a sum of connected embedded surfaces, $S = S_1 + \cdots + S_n$.

The next thing we wish to observe is that there is a nonzero degree map
from $S$ to the $n$-torus $\S_0$.   This map comes from the `almost
product' structure of $\hat M$ given in (\ref{coords}).    A simple deformation of the
$t$-lines in the vicinity of $K$ can be used to produce a continuous projection type map
$P: \hat M \to \S_0$
such that $K$ gets sent to a single point on $\S_0$ under $P$, and such that $P\circ j = {\rm id}$, where $j: \S_0 \to \hat M$ is inclusion.    Then $f = P \circ i : S \to \S_0$, where
$i : S \to \hat M$ is inclusion, is the desired nonzero degree map.  Indeed, $f$ induces
the map on homology $f_*: H_n(S) \to H_n(\S_0)$, and using that $S$ is homologous to $\S_0$, we compute, $f_*[S] =
P_*(i_*[S]) = P_*(j_*[\S_0]) = {\rm id}_*[\S_0] = \S_0 \ne 0$.

Thus, by linearity of $f_*$, at least one of the components of $S$, $S_1$, say,
admits a nonzero degree map to the $n$-torus.
By a result of
Schoen and Yau \cite{SY}, which does not require a spin assumption, $S_1$
does not admit a metric of positive scalar curvature.   (In fact it admits a metric
of nonnegative scalar curvature only if it is flat).
Moreover, we know that $S_1$ minimizes the brane action in its homology class
(otherwise there would exist a hypersurface homologous to $\S_0$
with brane action strictly less than that of $S$).  Thus, we can apply Theorem
\ref{warp} to conclude that a neighborhood $U$ of $S_1$ splits as a warped
product,
\beq
U = (-u_0, u_0) \times S_1 \qquad \hat g|_U = du^2 + e^{2u} h  \, ,
\eeq
where the induced metric $h$ on $S_1$ is flat.  But since $S_1$ in fact {\it globally}
maximizes the brane action in its homology class, by standard arguments this local
warped product structure can be extended to arbitrarily large $u$-intervals.  Hence
$K$ will eventually be contained in this constructed warp product region.  It now follows
that $\hat M$ has constant curvature $K_{\hat M} = -1$.  This in turn implies
that $M$ has constant curvature $K_{M} = -1$.  By previous remarks, we conclude
that $M$ is globally isometric to hyperbolic space.~\qed

\section{Positivity of mass}
The aim of this section is to give a proof of Theorem \ref{thm:posmass} on the
positivity of mass in the asymptotically hyperbolic setting.  We shall adopt here the definition of {\it asymptotically
hyperbolic}  given in Wang \cite{wang}:  
\begin{Def}\label{def:ah}
A  Riemannian manifold $(M^{n+1},g)$ is 
asymptotically hyperbolic provided it
is conformally compact, with
smooth conformal compactification $(\tilde M, \tilde
g)$, and with conformal boundary $\partial \tilde M = S^n$, such that
the metric $g$ on a deleted neighborhood $(0,T) \times S^n$ of $\partial
\tilde M = \{t = 0\}$ takes the form 
\beq\label{wmetric}
g = \sinh^{-2}(t) ( dt^2 +  h)  \,,
\eeq
where  $h = h(t,\cdot)$ is a
family of metrics on $S^n$, depending smoothly
on $t\in [0,T)$, of the form,
\beq
h = h_0 +t^{n+1} k + O(t^{n+2})  \,,
\eeq
where $h_0$ is the standard metric on $S^n$, and $k$ is a symmetric
$2$-tensor on $S^n$.
\end{Def}
We refer to $k$ as the mass aspect tensor; it is the leading order measure of the deviation of the metric $g$
from the hyperbolic metric.   Its trace with respect to $h_0$, 
${\rm tr}_{h_0}\,k$,
is called the mass aspect function.
%\footnote{the mass aspect function used in the work
%  of Wang differs from the one used here by a factor of $1/n$}, 
Up to a normalizing constant, the integral of the mass aspect function 
over the sphere defines the mass
(or energy)  of $(M,g)$, mass $= \int_{S^n} {\rm tr}_{h_0}\, k$.

For convenience,
we repeat here  the statement of our positivity of mass result.
\begin{thm} \label{thm:posmass2} 
Let $(M^{n+1},g)$, $2 \leq n \leq 6$, be an asymptotically hyperbolic
manifold with scalar curvature $S[g] \ge -n(n+1)$. Assume that the mass
aspect function 
$\tr_{h_0}\,k$ does not change sign, i.e. that it is either negative,
zero, or positive. Then, either the mass of $(M,g)$ is positive, or $(M,g)$
is isometric to hyperbolic space. 
\end{thm}

The proof, which makes use of the rigidity result Theorem \ref{rigid}, is carried out in the
next two subsections. In subsection \ref{deform} we obtain the deformation
result mentioned in the introduction, see Theorem \ref{defm} below.  This, together
with Theorem \ref{rigid}, implies that the mass aspect function cannot be negative,
see Proposition \ref{prop:aspect}.
In subsection \ref{rigidcase}, it is proved, using Theorem \ref{rigid}
again,  that if the mass aspect function vanishes then 
$(M,g)$ is isometric to hyperbolic space, see Theorem~\ref{thm:rigid2}.  These results together imply Theorem~\ref{thm:posmass2}.

\subsection{The deformation result}\label{deform}

Suppose $(M^{n+1},g)$ is  asymptotically hyperbolic in the sense of 
Definition \ref{def:ah}.   Then by  
making the  change of coordinate, $t = \arcsinh(\frac1{r})$, in \eqref{wmetric}
it follows that there is a relatively compact set $K$
such that   $M\setminus K = S^n \times [R,\infty)$, $R > 0$, 
and on $M \setminus K$, 
$g$ has  the form,
\beq\label{metric}
g = \frac1{1+r^2} \, dr^2 + r^2h \,, 
\eeq
where $h = h(\cdot, r)$ is an $r$-dependent family of metrics on $S^n$ of the form,
\beq
h = h_0+ \frac1{r^{n+1}}\, k + \sigma \,,
\eeq
where $h_0$ is the standard metric on $S^n$, $k$ is the mass aspect tensor and $\sigma = \sigma(\cdot,r)$ is an $r$-dependent family of metrics on
$S^n$ such that for integers $\ell, m \geq 0$,  one has,
\begin{equation} \label{eq:omder} 
|(r \partial_r)^{\ell}\d_x^m \sigma| \leq C/r^{n+2},  
\end{equation} 
for some constant $C$. 
For the proof of the deformation theorem and
the positive mass theorem, it is sufficient to assume condition
(\ref{eq:omder}) for $0 \leq \ell, m \leq 2$. 

Let $(M^{n+1},g)$ be asymptotically hyperbolic, with scalar curvature satisfying,
$S[g]\ge-n(n+1)$.
What we now prove is that if the mass aspect function of $(M,g)$ is pointwise negative then
$g$ can be deformed on an arbitrarily small 
neighborhood of infinity to the
hyperbolic metric, while preserving (after a
change of scale) the scalar curvature inequality $S \ge -n(n+1)$.
A more precise statement is given below. 
\begin{thm}\label{defm}
Let the metric $g$ be given as above. Suppose that the scalar
curvature of $g$, $S[g]$, satisfies $S[g]\geq -n(n+1)$. If the mass aspect
function 
$\tr_{h_0}\,k$ is pointwise negative, then for any sufficiently large $R_1 > R$ there
exists a metric $\hat g$ on $M$ such that,
\beq 
\hat g =  \begin{cases}  
g\,,  & \quad \mbox{$R \le r \le R_1$} \\\\
g_a = \frac{1}{1+\frac{r^2}{a}}\,dr^2+r^2h_0 \,, & \quad 
\mbox{$ 9\lambda R_1 \le r  < \infty $} \,,
\end{cases}
\eeq
where $a \in (0,1)$,  and such that, 
\beq\label{scalara}
S[\hat g] \ge -\frac{n(n+1)}{a} \,.
\eeq
(The constant $\lambda >1$ depends only on the mass aspect function; see 
section \ref{sec:preldef}.) 
\end{thm}
Theorems \ref{rigid} and \ref{defm} may be combined to give the following result.
\begin{prop}\label{prop:aspect}
Let $(M^{n+1},g)$, $2 \leq n \leq 6$, be an asymptotically hyperbolic
manifold with scalar curvature satisfying, $S[g] \ge -n(n+1)$. Then the mass
aspect function $\tr_{h_0}\,k$ cannot be everywhere 
pointwise negative.
\end{prop}

\begin{proof} 
Suppose to the contrary that the mass aspect function is strictly
negative.
Given any
$p\in M$, choose $R$ large enough so that $p$ is not in the end $S^n
\times[R, \infty)$.  Theorem \ref{defm}, together with a  
rescaling  of the metric, implies the existence of a metric $\tilde g$ on $M$
such that $(M,\tilde g)$ satisfies the hypotheses of Theorem \ref{rigid}.
Hence, $(M,\tilde g)$ is isometric to hyperbolic space $\bbH^{n+1}$.   But,
modulo the change of scale, by our construction, $\tilde g$ will differ from $g$
only at points on the end $S^n \times [R, \infty)$.  It follows that
$(M^{n+1},g)$ has constant negative curvature curvature in a neighborhood of
$p$. Since $p$ is arbitrary,   $(M^{n+1},g)$ must have globally  constant
negative curvature, which, by the asymptotics of $(M^{n+1},g)$, must equal
$-1$. Since $M$ is simply connected at infinity, we conclude that
$(M^{n+1},g)$ is isometric to hyperbolic space $\bbH^n$.  But this
contradicts the assumption that the mass aspect function is negative. 
\end{proof} 

\medskip
\noindent
{\it Proof of Theorem \ref{defm}:}
We now turn to the proof of the deformation result.
Introduce coordinates $x = (x^1, x^2, ..., x^n)$  on $S^n$.
We use the convention that for a function $f=f(x,r)$, $f'(x,r) = \partial_r
f(x,r)$, and $f''(x,r) = \partial_r^2 f(x,r)$.  

Let $\omega_{ij}$ and $k_{ij}$ be the components of $h_0$ and
$k$, respectively, with respect to the coordinates $(x^1, x^2,\dots , x^n)$.
Then $g$ in (\ref{metric}) takes the form
$$
g = \frac{dr^2}{1+r^2} + r^2\left(\omega_{ij} + \frac{\alpha_{ij}}{r^{n+1}}\right) dx^i dx^j
$$  
where 
$$
\alpha_{ij} = k_{ij}(x) + \frac{\beta_{ij}(x,r)}{r} \,.
$$ 
We are assuming (cf., \eqref{eq:omder}) $\alpha_{ij}$ satisfies the bounds, 
 \begin{align}\label{eq:alphabound} 
|(r\partial_r)^{\ell}\partial_x^m \alpha_{ij}| \le \Lambda   \,.
\end{align} 
for all integers $\ell, m$, $0 \le \ell,m \le 2$.

Let $\mu$ denote the the mass aspect function,
$\mu =  {\rm tr}_{h_0}k =   \omega^{ij} k_{ij}$;
by assumption,  $|\mu| > 0$. 
Let $\muup = \max_{x} |\mu(x)|$, $\mulow = \min_{x} |\mu(x)|$.
 We shall be making estimates of geometric quantities in terms of the above
defined constants. In particular, we shall use a generic constant
$$
C = C(n,R,\Lambda,\muup,\mulow)
$$
depending only on $n, R, \Lambda, \muup, \mulow$, which may change from line
to line. 
We shall further use the notation $O(1/r^k)$ for a quantity bounded by 
$C(n,R,\Lambda,\muup,\mulow)/r^k$.

\subsubsection{Preliminary definitions} \label{sec:preldef} 
Fix $R_1 > R$ to be specified later.
Set
$$
\lambda = \left ( \frac{\muup}{\mulow} \right )^{\frac{1}{n+1}}
$$
Let $a = a(n, \muup, \mulow, R_1) \in (0,1)$ be a number such that 
$$
\frac{n+1}{n} \sqrt{\frac{4}{3}} \frac{ \muup}{(4\lambda R_1)^{n+1}} 
< \frac{1}{a} - 1 < \frac{n+1}{n} 
\sqrt{\frac{3}{4}}\frac{ \mulow}{(3\lambda R_1)^{n+1}}
$$
 To show such an $a$ exists, it suffices to show that 
$$
\frac{n+1}{n} \sqrt{\frac{4}{3}} \frac{ \muup}{(4\lambda R_1)^{n+1}} 
< \frac{n+1}{n} \sqrt{\frac{3}{4}}\frac{ \mulow}{(3\lambda R_1)^{n+1}}
$$
or equivalently 
$$
\frac{3}{4} \left ( \frac{4}{3} \lambda \right )^{n+1} > \frac{\muup}{\mulow}  \,.
$$
By our choice of $\lambda$, this is equivalent to 
$$
\left ( \frac{4}{3} \right )^n > 1  
$$
which is obviously true. It follows from the definition that
 $a \nearrow 1$ as $R_1 \nearrow
\infty$.

It is
straightforward to show the existence of a smooth function 
$\psi: \Re \to \Re_+$ such that for any $R_1 > R$. 
\begin{subequations} \label{eq:psibounds} 
\begin{align}
\psi &= \left\{ \begin{array}{ll} 1, & r \leq 7\lambda R_1 \\ 
                                  0, & r \geq 8 \lambda R_1 
                 \end{array} \right.    \\
\psi'(r) &\leq 0 \text{ for all } r \\
|\psi'(r)| &\leq \frac{b}{r} \\
|\psi''(r)| &\leq \frac{c}{r^2}
\end{align}
\end{subequations} 
where  $b, c$ are positive constants. 
In the following, we consider a fixed function $\psi$ satisfying the
conditions (\ref{eq:psibounds}). 

In the deformation construction, we shall consider functions
$f: S^n \times [R, \infty) \to \Re$,  satisfying the
following conditions. 
\begin{subequations} \label{eq:fbounds} 
\begin{align} 
f(x,R) &= 1 , \text{ for } x \in S^n, \label{eq:fbound-1} \\ 
\intertext{ and for $(x,r) \in S^n \times [R, \infty)$ the conditions } 
\half  \leq  f &\leq 2 \label{eq:fbounda} \\ 
|f'(x,r)| &\leq  \frac{1}{r^2} \label{eq:fboundb} \\
|\Delta_{S^n} f| & \leq \frac{1}{r^n} \label{eq:fboundc} \\ 
|\nabla^{S^n} f|^2 & \leq \frac{1}{r^n} \label{eq:fboundd} \,.
\end{align} 
\end{subequations} 
In particular, the constant function $f \equiv 1$ satisfies conditions
(\ref{eq:fbounds}). 
In order to carry out the deformation from the metric $g$ to the metric $g_a$,
we shall consider metrics of the form,
\beq\label{metricdef}
g_{f,\psi} = \frac{1}{1+r^2 f} dr^2 + r^2 (\omega_{ij} + \frac{\psi
\alpha_{ij}}{r^{n+1}} ) dx^i dx^j    \,.
\eeq
Given $\psi$, the main objective is to construct an $f$ satisfying
\eqref{eq:fbounds} so that $g_{f,\psi}$ has the required properties.

\subsubsection{Scalar curvature formulas} 

We need formulas for the scalar curvature  $S[g_{f,\psi}]$ of the metric $g_{f,\psi}$.

\begin{lem} Let $f,\psi$ satisfy the assumptions (\ref{eq:fbounds}) and
(\ref{eq:psibounds}). Then the metric $g_{f,\psi}$ has scalar curvature,
\beq\label{scalargen}
S[g_{f,\psi}] =  -n(n+1)f-nrf'+ \frac1{r^n}(n|\mu|\psi' -  r|\mu|\psi'')f 
+ J  \,,
\eeq
where $J$ is a term bounded by $C(n, R,\Lambda,\muup,\mulow)/r^{n+2}$, 
and such that $J = 0$ for $r \ge 8\lambda R_1$.
\end{lem}
\begin{proof}
We describe our approach to carrying out this computation.
Setting,
\beq
h = 1 +r^2f \quad \mbox{and} \quad g_{ij} = r^2(\omega_{ij} + 
\frac{\psi\alpha_{ij}}{r^{n+1}}), 
\eeq
$g_{f,\psi}$ becomes,
\beq
g_{f,\psi} = \frac1{h} dr^2 + g_{ij}dx^idx^j \,.
\eeq

Applying the Gauss equation to an $r$-slice $\S = S^n \times \{r\}$ gives,
\beq\label{gauss}
S[g_{f,\psi}] = S_{\S}  + |B|^2 -H^2 +  2\ric(N,N)  \,,
\eeq
where $S_{\S}$, $B$ and $H$ are the  scalar curvature, second fundamental form
and mean curvature of $\S$, respectively, and $\ric(N,N)$ is the ambient Ricci
curvature in the direction of the unit normal $N=h^{1/2}\frac{\partial}{\partial r}$.
In terms of coordinates, $B$ and $H$ are given by, $b_{ij} = B(\d_i,\d_j) = \frac12 \sqrt{h}\,\d_r g_{ij}$,
and $H = g^{ij}b_{ij}$.  We then compute each term in (\ref{gauss}) in turn.

For the first three terms we obtain, making use of the bounds \eqref{eq:alphabound} and
\eqref{eq:psibounds}
\begin{align}
S_{\S} & =  \frac{n(n-1)}{r^2}+O(\frac{1}{r^{n+3}}) \label{induced} \\ \nonumber\\
H & =  h^{\frac12}[\frac{n}{r} +\frac{n+1}{2}\frac{|\mu|\psi}{r^{n+2}} - \frac{1}{2}\frac{|\mu|\psi^{\prime}}{r^{n+1}}+
O(\frac{1}{r^{n+3}})] \label{mean} \\ \nonumber\\
|B|^2 & =  B_i{}^jB_j{}^i = h
[\frac{n}{r^2}+(n+1)\frac{|\mu|\psi}{r^{n+3}}-\frac{|\mu|\psi'}{r^{n+2}}+O(\frac{1}{r^{n+4}})]
\label{square}\,.
\end{align}
Equations (\ref{mean}) and (\ref{square}) combine to give,
\begin{align}\label{combo}
|B|^2 - H^2 & =  h[-\frac{n(n-1)}{r^2}-(n-1)(n+1)\frac{|\mu|\psi}{r^{n+3}} \nonumber \\
& \hspace{1.2in} + (n-1)\frac{|\mu|\psi'}{r^{n+2}}+O(\frac{1}{r^{n+4}})]  \,.
\end{align}
Applying the Raychaudhuri  equation to the 
unit normal $N=h^{1/2}\frac{\partial}{\partial r}$ to the level sets of $r$, 
we have,
\beq\label{ric}
\ric(N,N)=-N(H)-|B|^2-\sqrt{h}\Delta_{\S}\frac{1}{\sqrt{h}} \,.
\eeq
Making use of equations (\ref{mean}) and (\ref{square}), we derive from (\ref{ric}),
\begin{align}\label{ric2}
\ric(N,N) & = -\frac{1}{2}h^{\prime}[\frac{n}{r} +
 \frac{n+1}{2}\frac{|\mu|\psi}{r^{n+2}} -\frac{1}{2}\frac{|\mu|\psi'}{r^{n+1}}+
O(\frac{1}{r^{n+3}})] 
\nonumber\\
&\quad + h[\frac{n(n+1)}{2}\frac{|\mu|\psi}{r^{n+3}} - n\frac{|\mu|\psi'}{r^{n+2}} +
\frac{1}{2}\frac{|\mu|\psi''}{r^{n+1}}+O(\frac{1}{r^{n+4}})]
\nonumber  \\ & \quad -\sqrt{h}\Delta_{\S}\frac{1}{\sqrt{h}}  \,.
\end{align}

Equations (\ref{induced}), (\ref{combo}) and (\ref{ric2}) then combine to give,
\begin{align}
S[g_{f,\psi}] &= \frac{n(n-1)}{r^2}-\frac{n(n-1)}{r^2}h-\frac{n}{r}\, h'  
%\nonumber\\
 - [\frac{n+1}{2}\frac{|\mu|\psi}{r^{n+2}}-\frac{1}{2}\frac{|\mu|\psi'}{r^{n+1}}
 + O(\frac{1}{r^{n+3}})]h^{\prime}
\nonumber\\
&\quad+ [(n+1)\frac{|\mu|\psi}{r^{n+3}} -(n+1)\frac{|\mu|\psi'}{r^{n+2}}+\frac{|\mu|\psi''}{r^{n+1}}+O(\frac{1}{r^{n+4}})]h
\nonumber\\
&\quad - 2 \sqrt{h} \triangle_{\S} \frac{1}{\sqrt{h}}  + O(\frac{1}{r^{n+3}})  \,.
\end{align}

Setting $h = 1 + r^2f$ in the above, and making use of the bounds
\eqref{eq:fbounds} and
\eqref{eq:psibounds}, one derives in a straight forward manner 
equation \eqref{scalargen}.  Moreover, it is clear from the computations that 
all `big O' terms vanish once $\psi$ vanishes.
\end{proof} 

The following Corollary gives the form of the scalar curvature which will be
 used in the deformation construction. 
\begin{cor}\label{cor:scalarineq}
 Let $f,\psi$ satisfy the assumptions (\ref{eq:fbounds}) and
(\ref{eq:psibounds}). Then, 
there is a nonnegative 
function 
$
A_1 : [R,\infty) \to \Re_+ , 
$
independent of $f$, 
such that 
\beq\label{A1}
A_1 \le C(n,R,\Lambda,\muup,\mulow) \quad \text{and} \quad A_1 = 0 
\quad \text{for $r \ge 9\lambda R_1$}   \, ,
\eeq
and such that the 
scalar curvature $S[g_{f,\psi}]$ of $g_{f,\psi}$ satisfies the
inequality 
\beq\label{scalarineq}
S[g_{f,\psi}] \geq - \frac{n}{r^n} \left [ 
(r^{n+1} + \frac{n+1}{n} |\mu| \psi - \frac{r}{n} |\mu| \psi') f -
\frac{A_1(r)}{r}\right ]'  \,.
\eeq
\end{cor} 
\begin{proof} 
Using the product rule, equation \eqref{scalargen} may be expressed as,
\begin{align}
S[g_{f,\psi}] & =  -\frac{n}{r^n} (r^{n+1}f)' + \frac1{r^n}[((n+1) |\mu|\psi -r|\mu|\psi')f]' 
\nonumber  \\ & \qquad\qquad  - \frac1{r^n}[(n+1) |\mu|\psi -r|\mu|\psi']f' +  O(\frac1{r^{n+2}}) \nonumber \\ \intertext{which, by the bounds  (\ref{eq:fbounds}) and
(\ref{eq:psibounds}) simplifies to,} 
S[g_{f,\psi}]  & =  -\frac{n}{r^n} [(r^{n+1} - \frac{n+1}{n}|\mu|\psi + \frac{r}{n}|\mu|\psi')f]' 
+   O(\frac1{r^{n+2}})  \,.
\end{align} 
It follows that there exists a smooth function $A: [R, \infty) \to \bbR_+$,
satisfying,   $A \le C(n,R,\Lambda,\muup,\mulow)$ and $A = 0$ 
for $r \ge 9\lambda R_1$, such that,
\beq\label{scalarineq2}
S[g_{f,\psi}] \ge    -\frac{n}{r^n} [(r^{n+1} - \frac{n+1}{n}|\mu|\psi + \frac{r}{n}|\mu|\psi')f]' - 
\frac{A(r)}{r^{n+2}}.
\eeq
Now define $A_1 : [R,\infty) \to \Re_+$ by,
\beq
A_1(r) =  \frac{r}{n} \int_r^{9\lambda R_1} \frac{A(t)}{t^2} dt  \,.
\eeq
One easily checks that the properties \eqref{A1} hold. Moreover, since,
\beq
\left(\frac{A_1(r)}{r}\right)' = -\frac1{n}\frac{A(r)}{r^2}\,,
\eeq
inequality \eqref{scalarineq} follows from \eqref{scalarineq2}. Finally, it
is clear from the construction that $A$ and hence also $A_1$ may be chosen
to be independent of $f$. 
\end{proof}

\subsubsection{Defining $\eta$: rounding the corner}\label{sec:etadef}  

Let $\eta_1(x,r)$, $\eta_2(r)$ be given by 
\begin{subequations} \label{eq:eta12def}
\begin{align}
\eta_1 (x,r) &= r^{n+1} + \frac{n+1}{n} |\mu| \psi - \frac{r}{n} |\mu| \psi' -
\frac{A_1(r)}{r}, \\
\eta_2(r) &= \frac{r^{n+1}}{a},
\end{align} 
\end{subequations} 
where the function $A_1$ appearing in $\eta_1$ is the function $A_1$
defined in Corollary \ref{cor:scalarineq}. 
By \eqref{scalarineq},  $S[g_{1,\psi}] \ge - \frac{n}{r^n} \eta_1'$. Further,
$g_a$ has constant curvature $-1/a$, and hence 
$S[g_a] =  - \frac{n}{r^n} \eta_2'$.

In the rest of the argument we shall be choosing $R_1$ sufficiently large, so
that the required conditions are satisfied. We shall successively increase
$R_1$ as required. 

\begin{lem} \label{lem:eta12est}
There is an $R_1 > R$, $R_1 = R_1(n,R,\Lambda,\muup,\mulow)$ 
such that the following inequalities hold.
\begin{subequations}\label{eq:etaineq}
\begin{align} 
\eta_1'(x,r) &< \eta_2'(r) \,, && R_1 \leq r \leq 7\lambda R_1
\label{eq:etaineq3}  \\
\eta_1(x,r) - \eta_2(r) &> \mulow/20 \,,
&& R_1 \le r \leq 3\lambda R_1   \label{eq:etaineq1} \\ 
\eta_2(r) - \eta_1(x,r) & > \muup/10 \,, 
&& 4\lambda R_1  \leq r \leq 7 \lambda R_1 \label{eq:etaineq2}  \,.
\end{align} 
\end{subequations} 
%(See Figure 1.)
\end{lem} 
\begin{figure}[Ht]
\psfrag{META1}{$\eta_1$}
\psfrag{META2}{$\eta_2$}
\psfrag{META}{$\eta$}
\psfrag{MR}{$R_1$}
\psfrag{M3R}{$3\lambda R_1$}
\psfrag{M4R}{$4\lambda R_1$}
\psfrag{M7R}{$7\lambda R_1$}
\centering 
\includegraphics[width=2.5in]{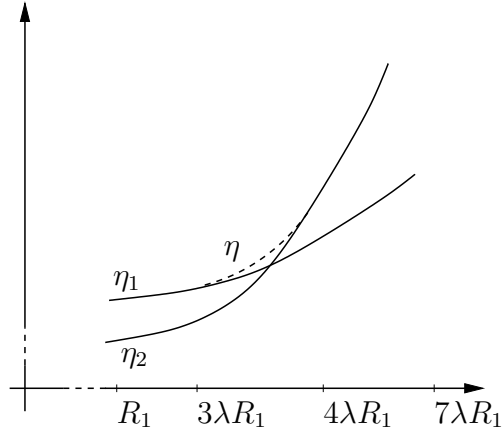} 
\caption{Construction of $\eta$}
\label{fig2}
\end{figure}
\begin{proof} 
We need to consider only the interval $R_1 \leq r \leq 7 \lambda R_1$. There,
$\psi \equiv 1$, and $\psi' \equiv 0$. 
We start by proving (\ref{eq:etaineq3}). By construction, we have for $r \geq
R_1$, using $|\mu| > 0$ and the properties of $a$, after some manipulations, 
\begin{align*} 
\eta_2' - \eta_1' &> \frac{1}{r} \left [ 
(\frac{1}{a} - 1 ) (n+1) r^{n+1} - \frac{A(r)}{nr} \right ] \\
& > \frac{1}{r} \left [ \sqrt{\frac{4}{3}} \frac{(n+1)^2}{n}
  \frac{\muup}{(4\lambda)^{n+1}} 
-   \frac{A(r)}{nr} \right ] 
\end{align*} 
Since $A(r) \leq C(n,\Lambda, R, \muup,\mulow)$, it follows that 
(\ref{eq:etaineq3}) holds for $R_1$ sufficiently large. 

Next we prove (\ref{eq:etaineq1}). Since $(\eta_1 - \eta_2)' < 0$ for $r \geq
R_1$, it is sufficient to prove the inequality at $r = 3\lambda R_1$. We have 
\begin{align*} 
\eta_1(x,3\lambda R_1) - \eta_2(3\lambda R_1) &\geq 
\left [ \left (1 - \frac{1}{a} \right ) + \frac{n+1}{n} \frac{\mulow}{(3\lambda R_1)^{n+1}}
  \right ] (3\lambda R_1)^{n+1} \\
&\quad - \frac{A_1(3\lambda R_1)}{3\lambda R_1}  \,,\\ 
\intertext{which by using the properties of $a$ and simplifying gives,}
&> \left [ 1- \sqrt{\frac{3}{4}} \right ] \frac{n+1}{n} \mulow - \frac{A_1(3\lambda R_1)}{3\lambda R_1}   \,.
\end{align*} 
One checks that $1 - \sqrt{3/4} > 1/10$. 
Thus, in view of the fact that $A_1 \leq
C(n,R,\Lambda, \muup, \mulow)$, by possibly increasing $R_1$ we see
that (\ref{eq:etaineq1}) can be made to hold.

We proceed in a similar fashion to prove  (\ref{eq:etaineq2}). 
We have 
\begin{align*} 
\eta_2(4\lambda R_1) - \eta_1(x,4\lambda R_1) &\geq 
\left [ (\frac{1}{a} - 1) - \frac{n+1}{n} \frac{\muup}{(4\lambda R_1)^{n+1}} 
\right ] (4\lambda R_1)^{n+1} \\
&\quad \geq \left [ \sqrt{\frac{4}{3}} -1 \right ] \frac{n+1}{n} \muup
\end{align*} 
We note that  $\sqrt{4/3}-1 > 1/10$. 
 This completes the proof of Lemma \ref{lem:eta12est}. 
\end{proof} 

We shall define $\eta = \eta(x,r)$ to be a suitably increasing function of the
variable $r$ that
smoothly transitions from $\eta_1$  to $\eta_2$ (see Figure 1).  In order to 
make meaningful estimates, its construction shall be made fairly
explicit.  Its construction depends on two auxiliary functions $\a$ and
$\b$, which we now introduce.

Let $\alpha: [0, \infty) \to \Re$ be a function satisfying, 
\begin{enumerate} 
\item 
$\alpha(r) = 0$ for $r \leq R_1$, $2R_1 \leq r \leq 5\lambda R_1$, and $r
\geq 6\lambda R_1$. 
\item $\a > 0$ for $R_1 < r < 2R_1$, $\a < 0$ for $5\lambda R_1 < r < 6\lambda R_1$,  
and
\begin{equation}\label{eq:intzalp}
\int_{R_1}^{6\lambda R_1} \alpha(t) dt = \int_{-\infty}^\infty \alpha(t) dt =
0  \,.
\end{equation} 
\end{enumerate} 
Consider,
$$
\gamma(x,r) = \frac{\alpha(r)}{\eta_1(x,r) - \eta_2(r)}   \,.
$$
It follows from Lemma \ref{lem:eta12est} and the properties of $\alpha$ that
$\gamma$ is nonnegative and bounded, 
$$
0 \leq \gamma \leq C(n,\Lambda, \mulow,\muup)  \,.
$$
Next, define $m(x)$ by the condition 
\beq\label{mdef}
1 + m(x) \int_{R_1}^{6\lambda R_1} 
%\frac{\alpha(t)}{\eta_1(x,t) -\eta_2(x,t)} 
\g(x,t)\, dt = 0  \,,
\eeq
and let 
$$
\beta(x,r) = 1 + m(x) \int_{R_1}^{r}
% \frac{\alpha(t)}{\eta_1(x,t) - \eta_2(x,t)}   
\g(x,t)\, dt \,.
$$
Then, $\b$ satisfies,
\begin{subequations}\label{eq:betaprop} 
\begin{align} 
\beta &= 1, && \text{ for $r \leq R_1$} \label{eq:betaprop1} \\ 
\beta &= 0, && \text{for $r \geq 6 \lambda R_1$} \label{eq:betaprop2} \\ 
0 \leq \beta &\leq 1 , && \text{ for all $r$} \label{eq:betaprop3}  \,.
\end{align} 
\end{subequations} 
Conditions (\ref{eq:betaprop1}), (\ref{eq:betaprop2}) are clear. For
(\ref{eq:betaprop3}), we consider,  
$\beta' = m\g$.  
As observed above, $\g \ge 0$,  and hence from \eqref{mdef},
$m \le 0$, so that $m\g \le 0$.
Hence $\beta$ is decreasing, which implies that $0 \leq \beta
\leq 1$. 

Now we are ready to define $\eta$. Let 
$$
\eta(x,r) = \eta_1(x,R_1) + \int_{R_1}^r  [ \beta(x,t) \eta_1'(x,t)  + (1-\beta(x,t))  
\eta_2'(t) ] dt    \,.
$$
\begin{lem} 
$\eta$ defined as above satisfies the conditions 
\begin{subequations} 
\begin{align} 
\eta(x,r) &= \left \{ \begin{array}{ll} 
\eta_1(x,r)
& r \leq R_1  \label{eq:etapropa} \\ \\
\eta_2(r) = \frac{r^{n+1}}{a} , & r \geq 6 \lambda R_1 \end{array} \right. \\ \nonumber \\
\eta' (x,r) &\leq \eta_2'(r) =  \frac{(n+1)r^n}{a} , \quad R_1 \leq r  < \infty     \,.
%\leq 8 \lambda R_1
\label{eq:etapropb} 
\end{align} 
\end{subequations}
\end{lem} 
\begin{proof} 
For $r \leq R_1$, we have $\beta = 1$, so 
$$
\eta(x,r) = \eta_1(x,R_1) + \int_{R_1}^r \eta_1'(x,t) dt = \eta_1(x,r)  \,.
$$
In the following calculations, which only involve derivatives and integrals
with respect to $r$, we suppress reference to $x$ in order to avoid
clutter. 
For $r \geq 6 \lambda R_1$, we have $\beta = 0$, which gives 
\begin{align} 
\eta(r) &= \eta_1(R_1) 
+ \int_{R_1}^{6\lambda R_1} [ \beta  \eta_1' +
(1-\beta) \eta_2' ] dt 
+ \int_{6\lambda R_1}^r \eta_2' dt \\
&= \eta_2(r) - \eta_2(6\lambda R_1) + \eta_1(R_1) 
+ \int_{R_1}^{6\lambda R_1} (\beta \eta_1' + (1-\beta) \eta_2') dt  \,.
\end{align}   
A partial integration gives 
\begin{align} 
\int_{R_1}^{6\lambda R_1} & ( \beta \eta_1' + (1-\beta) \eta_2' ) dt 
\\
&= \eta_2(6\lambda R_1) - \eta_2(R_1) + \beta (\eta_1 - \eta_2)
\bigg{|}_{r=R_1}^{r=6\lambda R_1} - \int_{R_1}^{6\lambda R_1} \beta' (\eta_1 -
\eta_2) dt \\
\intertext{use the properties of $\beta$ and $\gamma$} 
%section \ref{sec:betadef},} 
&= \eta_2(6\lambda R_1) - \eta_2(R_1) - (\eta_1(R_1) - \eta_2(R_1) ) 
- m \int_{R_1}^{6\lambda R_1} \alpha dt \\
\intertext{use (\ref{eq:intzalp})}
%see section \ref{sec:etadef}} 
&= \eta_2(6\lambda R_1) - \eta_1(R_1)  \,.
\end{align} 
Substituting this into the formula for $\eta(x,r)$, we obtain,
$$
\eta(x,r) = \eta_2(r), \quad \text{ for $r \geq 6\lambda R_1$} 
$$
We have now established (\ref{eq:etapropa}). 
Next we prove (\ref{eq:etapropb}). We have 
\begin{align} 
\eta' &= \beta \eta_1' + (1-\beta) \eta_2' \\
\intertext{and hence, by \eqref{eq:etaineq},} 
%and \eqref{eq:etapropa},}
&\eta' \leq \beta \eta_2' + (1-\beta) \eta_2' \\ 
&= \eta_2' = \frac{(n+1)r^n}{a}   \,.
\end{align} 
\end{proof} 
We shall now make use of the function $\eta$ defined above to define a
function $f$, which shall be shown to satisfy 
the conditions (\ref{eq:fbounds}). This fact allows us to apply the result of
Corollary \ref{cor:scalarineq} to estimate the scalar curvature of the
deformed metric $g_{f,\psi}$ defined in terms of this $f$. 

\subsubsection{Defining $f$} 
We define $f$ by, 
\begin{equation}\label{eq:f-eta}
f = \frac{\eta + \frac{A_1}{r}}{\eta_1 + \frac{A_1}{r}}   \,,
\end{equation} 
which implies,
\begin{equation}\label{eq:rAeta} 
(r^{n+1} + \frac{n+1}{n} |\mu| \psi - \frac{r}{n} |\mu| \psi' ) f - \frac{A_1}{r}
= \eta  \,.
\end{equation} 
Here, $A_1 = A_1(r)$ is the function determined in Corollary
\ref{cor:scalarineq}. It is crucial to note here that $A_1$ is {\em
  independent} of the particular $f$, as long as it satisfies the conditions
(\ref{eq:fbounds}). Our task is now to show that $f$ defined as
above does satisfy these conditions as long as $R_1$ is chosen sufficiently
large. This will be demonstrated in Lemma \ref{lem:fOK} below. It then
follows from equation (\ref{eq:rAeta}) above, Corollary~\ref{cor:scalarineq}
and  \eqref{eq:etapropb} that 
$S[g_{f,\psi}]$ satisfies,
$$
S[g_{f,\psi}] \geq -\frac{n}{r^n}\eta'  \ge  -\frac{n(n+1)}{a} \,.
$$
In addition we have, $g_{f,\psi} = g$ on $[R,R_1]$ and $g_{f,\psi} =g_a$
on $[9\lambda R_1,\infty)$. 
Thus, 
subject to the following lemma, Theorem \ref{defm} has been proven.

\begin{lem} \label{lem:fOK} 
Let $\eta$ be as in section \ref{sec:etadef}, and let $f$ be given in terms of $\eta$ by
(\ref{eq:f-eta}). Then, there is an $R_1=R_1(n,R,\Lambda,\muup,\mulow)$ 
sufficiently large, so that the  
inequalities (\ref{eq:fbounds}) are valid.   
\end{lem} 
\begin{proof} 
The condition (\ref{eq:fbound-1}) is clear from the construction. Since
(\ref{eq:fboundb}) implies (\ref{eq:fbounda}) we only need to verify 
that $|f'| \leq 1/r^2$ and $|\partial_x f| \leq 1/r^n$, $|\partial_x^2 f |
\leq 1/r^n$. 
We begin by showing there is an $R_1$ sufficiently large, and not smaller
  than the previously made choices of $R_1$, 
so that $|f'| \leq 1/r^2$. 
We have, 
\begin{equation}\label{eq:f-eta1}
f = 1 + \frac{(\eta - \eta_1)}{\eta_1 + A_1/r}   \,.
\end{equation} 
This gives 
$$
f' = \frac{(\eta - \eta_1)'}{\eta_1 + A_1/r} 
- \frac{(\eta - \eta_1)}{(\eta_1 + A_1/r)^2}[\eta_1' + (A_1/r)']  \,.
$$
Since 
\begin{align*} 
\eta(x,r) &= \eta_1(x,R_1) + \int_{R_1}^r (\beta \eta_1' + (1-\beta) \eta_2')  
dt  \,, \\
\eta_1(x,r) &= \eta_1(x,R_1) + \int_{R_1}^r \eta_1' dt \, \\
\intertext{we have,}
\eta - \eta_1 &= \int_{R_1}^r (1-\beta) (\eta_2 - \eta_1)'  \,.
\end{align*} 
so 
$$
(\eta - \eta_1)' = (1-\beta) (\eta_2 - \eta_1)' 
$$
Hence, 
\begin{align} 
|(\eta - \eta_1)'| 
\leq (\eta_2 - \eta_1)' &= \left (\frac{1}{a} - 1 \right ) r^n (n+1) + O(1/r)   
\\
&\leq \frac{C}{r}  \,.
\end{align} 
Here and below $C = C(n,R,\Lambda,\muup,\mulow)$ is a generic constant. 
Recall that by construction we have 
$$
(A_1/r)' \leq \frac{C}{r^2}   \,.
$$
Combining the above inequalities, we have 
\begin{equation} \label{eq:*} 
\left | \frac{(\eta -\eta_1)'}{\eta_1 + A_1/r} \right | 
\leq \frac{C}{r^{n+2}}  \,.
\end{equation} 
Now since 
$$
|(\eta - \eta_1)'| \leq \frac{C}{r}
$$
we have 
\begin{equation}\label{eq:eta-eta1-bound}
|\eta - \eta_1| \leq C
\end{equation} 
for $r \in [R_1 , 9 \lambda R_1]$. This together with 
$$
|\eta_1'| \leq C r^n 
$$
implies that 
\beq\label{eq:**}
\left | \frac{(\eta - \eta_1)}{(\eta_1 + A_1/r)^2} [\eta_1' +
(A_1/r)'] \right | \leq \frac{C}{r^{n+2}}   \,.
\eeq
Equations (\ref{eq:*}) and (\ref{eq:**}) imply 
%\mnote{missing ref here}
$$
|f'| \leq \frac{C}{r^{n+2}}
$$
for $r \in [R_1, 9 \lambda R_1]$. By choosing $R_1$ large enough, we have 
$$
|f'| \leq \frac{1}{r^2}
$$ 
which gives (\ref{eq:fboundb}).

Next we demonstrate that by, if necessary, further increasing $R_1$, we can
ensure that the condition 
$$
|\partial_x^k f| \leq \frac{1}{r^n} , \quad k = 1,2,
$$
holds, 
where $\d_x$ denotes partial differentiation with respect to any one of the
coordinates $x^i$.

Recalling (\ref{eq:f-eta1}), we have  
$$
\partial_x f = \frac{\partial_x(\eta - \eta_1)}{\eta_1 + A_1/r} 
- \frac{\eta - \eta_1 }{(\eta_1 + A_1/r)^2}(\partial_x \eta_1) 
$$
We estimate the second term first. By (\ref{eq:eta-eta1-bound}), $|\eta - \eta_1|
\leq C$, so 
$$
\left | \frac{\eta - \eta_1}{(\eta_1 + A_1/r)^2} \partial_x \eta_1
\right | \leq \frac{1}{r^{2n+2}} |\partial_x \eta_1| 
$$
Recalling
$$
\eta_1 = r^{n+1} + \frac{n+1}{n} |\mu| \psi - \frac{r}{n} |\mu| \psi' -
\frac{A_1}{r} 
$$
we find  $|\partial_x \eta_1| \leq C$. Hence
the modulus of the second term in $\partial_x f$ is bounded by $C/r^{2n+2}$. 

Now consider the first term in $\partial_x f$. Recall 
\begin{align*} 
\eta - \eta_1 &= \int_{R_1}^r (1-\beta) (\eta_2 - \eta_1) ' dt \\
&= (1-\beta) (\eta_2 - \eta_1) \bigg{|}_{R_1}^r + \int_{R_1}^r (\eta_2 -
\eta_1) \beta ' dt \\
&= [1- \beta] (\eta_2 -  \eta_1) + \int_{R_1}^r (\eta_2 - \eta_1) \beta' dt 
\end{align*} 
so 
\begin{multline*}
\partial_x (\eta - \eta_1) = - (\partial_x \beta) (\eta_2 - \eta_1) - (1-\beta)
\partial_x \eta_1 \\
+ \int_{R_1}^r \left ( [ \partial_x (\eta_2 - \eta_1) ] \beta' + (\eta_2 -
\eta_1) \partial_x \beta' \right ) dt   \,.
\end{multline*} 
Now,
\begin{align*} 
\partial_x \beta &= 
(\partial_x m) \left ( \int_{R_1}^r \frac{\alpha}{\eta_1 -
\eta_2} \right ) 
+ m \int_{R_1}^r \left ( \partial_x \frac{\alpha}{\eta_1 - \eta_2} \right ) \,.
\end{align*} 
Recall, 
$$
1+ m \int_{R_1}^{6\lambda R_1} \frac{\alpha}{\eta_1 - \eta_2} = 0
$$
so 
$$
m = - \frac{1}{\int_{R_1}^{6\lambda R_1} \frac{\alpha}{\eta_1 - \eta_2}} = O(1/R_1)  \,.
$$
Hence,  
$$
\partial_x m = - \frac{\int_{R_1}^{6\lambda R_1} \partial_x \left (
\frac{\alpha}{\eta_1 - \eta_2} \right )}
{\left ( \int_{R_1}^{6\lambda R_1} \frac{\alpha}{\eta_1 - \eta_2} \right )^2 }  \,.
$$
One can see that, 
$$
\partial_x \left ( \frac{\alpha}{\eta_1 - \eta_2} \right ) = O(1) 
$$
so 
$$
\partial_x m = O(\frac{1}{R_1})  \,.
$$
Hence we have, 
using $m = O(1/R_1)$,  
$$
\partial_x \beta = O(1)
$$
and therefore 
$$
\partial_x \beta (\eta_2 - \eta_1) = O(1) \,.
$$
Next we consider the second term in $\partial_x f$. 
We have 
$$
|\partial_x \eta_1 | \leq C \,,
$$
so 
$$
(1-\beta) \partial_x \eta_1 = O(1)  \,.
$$
Since $m = O(1/R_1)$, we have, 
$$
\beta' = m \frac{\alpha}{\eta_1 - \eta_2} = O(1/R_1) 
$$
but, 
$$
\partial_x (\eta_2 - \eta_1) = O(1) 
$$
so, 
$$
\partial_x (\eta_2 - \eta_1) \beta' = O(1/R_1)  \,.
$$
Similarly, 
$$
\partial_x \beta ' = O(1/R_1)
$$
and, 
$$
\eta_2 - \eta_1 = O(1)   \,.
$$
Hence, 
$$
[\partial_x ( \eta_2 - \eta_1) ] \beta' + (\eta_2 - \eta_1) \partial_x \beta'
= O(1/R_1)   \,,
$$
and therefore, 
$$
\int_{R_1}^r [ (\partial_x ( \eta_2 - \eta_1) ) \beta' + (\eta_2 - \eta_1)
\partial_x \beta' ] = O(1)   \,.
$$
Adding the three terms, we get, 
$$
\partial_x (\eta - \eta_1) = O(1)  \,.
$$
This implies that the first term in $\partial_x f$ is 
$$
\frac{O(1)}{\eta_1 + A_1/r} = O(\frac{1}{r^{n+1}})  \,,
$$
i.e.,
$$ 
|\partial_x f| \leq \frac{C}{r^{n+1}}  \,.
$$
Similar arguments give,  
$$
|\partial_x^2 f| \leq \frac{C}{r^{n+1}}   \,.
$$
By choosing $R_1$ large enough we obtain,  
$$
|\partial_x^k f| \leq \frac{1}{r^n} , \quad k=1,2   \,.
$$
Lemma \ref{lem:fOK} follows. 
\end{proof} 
As discussed above, now that Lemma \ref{lem:fOK} is established, we have
completed the proof of Theorem
\ref{defm}.

\subsection{The case of vanishing mass aspect function}\label{rigidcase}

In this subsection, we prove the following.

\begin{thm} \label{thm:rigid2} 
Let $(M^{n+1},g)$, $2 \leq n \leq 6$, be an asymptotically hyperbolic
manifold with scalar curvature satisfying, $S[g] \ge -n(n+1)$. If  the mass
aspect function 
$\tr_{h_0}\,k$ vanishes identically, then $(M,g)$ 
is isometric to hyperbolic space. 
\end{thm}
We note  that while Theorem \ref{thm:rigid2} generalizes Theorem \ref{rigid},
its proof relies on it.  We note also that our positivity of mass 
result, Theorem \ref{thm:posmass2}, follows immediately from Proposition \ref{prop:aspect} and Theorem
 \ref{thm:rigid2}.

For notational convenience we set $d = n+1$. Further, 
let capital latin indices run from $1,\dots,d-1$, let lowercase latin indices
run from $1, \dots, d$, and let $y^A$ be coordinates on $S^{d-1}$. Further,
let  $(x^i)
= (t, y^A)$  be coordinates on $(0,T) \times S^{d-1}$, and, as usual,   
let $h_0$ be the standard metric on $S^{d-1}$.

\subsubsection{Conformal gauge} 

Consider a conformally compact $d$-dimensional manifold $(M,g)$ where $M$ is
the interior of a manifold with boundary $\tM = M \cup \partial M$, 
and suppose  $g$ is of the form $g = \rho^{-2} \tme$, 
with $\rho$ a defining function for $\partial M$ for a metric $\tme$ which is
smooth on $\tM$. 

Let $\theta$ be a positive function on $\tM$. Letting $\tme \to \theta^2 \tme$
and $\rho \to \theta \rho$ leaves $g$ unchanged. Such a transformation can
therefore be viewed as a change of conformal gauge.

Let $\tme$ be a metric on $\tilde M$ which in a neighborhood 
of $\partial M$
can be written in the form 
\begin{equation}\label{eq:tmeform}
\tme = dt^2 + h_0 + t^d \gamma 
\end{equation} 
where $\gamma = \gamma_{ij} dx^i dx^j$ is a smooth tensor field on $(0, T)
\times S^{d-1}$ for some $T > 0$, such that 
the restriction of $\gamma$ to $\partial M$ is a smooth tensor on $S^{d-1}$, 
i.e. $\gamma\big{|}_{\partial M} = \gamma(0,y)_{AB} dy^A dy^B$. The following
lemma shows that after a change of conformal gauge we may assume that $\tme$ is
in Gauss coordinates based on $\partial M$. 
\begin{lem} \label{lem:gauge} 
Consider the conformally compact metric $g = \sinh^{-2}(t) \tme$. 
There is a conformal gauge change so that $g$ takes the form $g =
\sinh^{-2}(\that) 
%\hat
{\htme}$, where $\that(p) = d_{\htme}(p, \partial M)$, and 
$\htme$ is of the form 
$$
\htme = d\that^2 + h_0 + \that^d \hgamma 
$$
where $\hgamma = \hgamma(\hat t,y)_{AB} dy^A dy^B$ is a $\that$-dependent tensor field
on $S^{d-1}$, such that $\hgamma \big{|}_{\partial M} = \gamma
\big{|}_{\partial M}$. In particular, $g$ is asymptotically hyperbolic in the sense
of Definition \ref{def:ah}, with 
mass aspect tensor $k = \gamma \big{|}_{\partial M}$.
\end{lem} 
\begin{proof} 
Let $\rho = \sinh(t)$. 
Arguing as in \cite[Section 5]{AD}, 
we shall find a function $\theta$ such that $\hat \rho :
= \theta \rho = \sinh(\that)$, where $\that = d_{\htme}(p,\partial M)$, 
is the distance to the boundary
in the metric $\htme$. This is equivalent to the condition that 
$\hat f := \arcsinh(\hat \rho) = \that$, with $|d\hat f|_{\htme} = 1$. A
calculation as in the proof of \cite[Lemma 5.3]{AD} shows that this condition
is equivalent to the equation 
\begin{equation}\label{eq:diff} 
\rho \tme (d\theta,d\theta) + 2 \theta \tme (d\theta, d\rho) = \theta^4 \rho
+ \theta^2 a
\end{equation}
where $a = \rho^{-1} (1-\tme(d\rho, d\rho))$. 

For $\tme$ of the form (\ref{eq:tmeform}), we have $a = -\rho +
O(t^d)$. Equation (\ref{eq:diff}) is a system of first order partial differential
equations, with characteristics transversal to $\partial M = \{ t = 0\}$, and
satisfies the conditions for existence of solutions with initial condition
$\theta = 1$ at $\partial M$, see \cite[volume 5, pp. 39-40]{spivak}. Hence
there is a small neighborhood $U$ of $\partial M$, and a solution $\theta$ to
(\ref{eq:diff}) on $U$. 

We shall need the following fact. 
\begin{claim} $\theta = 1 + t^{d+1}w$, where $w$ is smooth up to $\partial
  M$
\end{claim}
The proof of the claim is straightforward and is left to the reader. 
Now we have that $\that = \arcsinh(\theta\sinh(t)) = t [ 1 + O(t^{d+1})]$, and
hence 
$$
t = \that [ 1 + O(\that^{d+1})]
$$
where the $O(t^{d+1})$ and $O(\that^{d+1})$ terms are smooth functions of
$(t,y)$ and $(\that, y)$, respectively. 
It is  straightforward to verify that  
$$
\sinh^{-2}(t) = \sinh^{-2}(\that)[1 + \that^{d+1}]
$$
and $g = \sinh^{-2}(\that) \htme$, with 
$$
\htme = d\that^2 + h_0 + \that^d \hgamma 
$$
where $\hgamma$
has the property that 
$\hgamma(0,y) = \gamma(0,y)$.

By construction, $\that$ is the distance to $\partial M$, and hence the above
is the form of $\htme$ in Gauss coordinates, based on $\partial M$. It
follows that $\hgamma$ is a $\that$-dependent tensor on $S^{d-1}$. 
\end{proof}

\subsubsection{Conformal deformation} 

Assume $(M,g)$ is asymptotically hyperbolic in the sense of Definition
\ref{def:ah}.  
Then, in slightly different notation, $(M,g)$
has a conformal compactification 
$(\tM, \tme)$ with conformal boundary $\partial \tM$ 
the round sphere, and such that near $\partial \tM$, $g$ has the form, 
\beq\label{rhometric}
g = \rho^{-2} \tme  
\eeq
where $\rho = \sinh(t)$ and 
\beq\label{rhometric2}
\tilde g = dt^2 + h_0 + t^d \gamma \,,
\eeq
where $h_0$ is the standard metric on
$S^{d-1}$ and $\gamma = \g(t,\cdot)$ is a $t$-dependent family of metrics
on $S^{d-1}$
smooth up to $\partial \tM$. Note that the mass aspect tensor is given by, $k = 
\gamma_{AB} \big{|}_{\partial \tM}$.

Let $\hpara = h_0 + t^d \gamma$ be the metric
induced on the level sets of $t$, and let $\Spara$
denote the  scalar curvature defined with
respect to $\hpara$. Further, let $K_{ij} = \half \partial_t h_{ij}$. The
only nonvanishing components of $K$ are 
$K_{AB} = \half d t^{d-1} \gamma_{AB} + O(t^d)$. 

Let $\tnabla$ denote the covariant derivative defined with respect to $\tme$
and let $\tilde S$ denote the scalar curvature of $\tme$.  The formula 
for the scalar curvature of conformally related metrics gives,
$$
S = -d(d-1) \tnabla_l \rho \tnabla^l \rho + (2d-2) \rho \tnabla^l \tnabla_l 
\rho + \rho^2 \tS   \,.
$$
\begin{claim}  $S$ has the asymptotic form,
\beq\label{eq:scal}
S = -d(d-1) +  O(t^{d+1})    \,.
\eeq
\end{claim}
Indeed, by Taylor's theorem, we have  $\Spara = S[h_0] + O(t^d)$, and hence,
\begin{align*} 
\tS &= \Spara - 2 h_0^{AB} \partial_t K_{AB} - (h_0^{AB} K_{AB})^2 + 3 K_{AB}
K^{AB}  \\
&=(d-1)(d-2) - d(d-1) t^{d-2} h_0^{AB} \gamma_{AB} + O(t^{d-1})    \,.
\end{align*}

Further,  using $\rho = \sinh t$,
\begin{align*}
\tnabla^l \tnabla_l \rho &= \sinh(t) - \tme^{ij} \tGamma_{ij}^t \cosh(t) \\
&= \sinh(t) + h_0^{AB} K_{AB} + O(t^d) \\ 
&= \sinh(t) + \frac{d}{2} t^{d-1} h_0^{AB} \gamma_{AB} + O(t^d)  \,.
\end{align*} 
Finally, we note that $\tnabla^l \rho \tnabla_l \rho = \cosh^2(t)$. Putting
this together, one finds after a few manipulations 
that the terms involving the mass aspect function, $\mu = h_0^{AB} ( \gamma_{AB} \big{|}_{\partial \tM})$,
in $S[g]$, at order $t^d$, cancel.  Equation \ref{eq:scal} follows.

\smallskip
By standard results \cite{aviles:mcowen}, 
there is a unique positive
solution 
$u$ such that $\lim_{x \to \infty}
u(x) = 1$,
to the Yamabe equation  for prescribed scalar curvature $-d(d-1)$ 
in dimension $d$,
\begin{equation}\label{eq:yam}
- \frac{4(d-1)}{d-2} \Delta u + S u + d(d-1) u^{\frac{d+2}{d-2}} = 0  \,.
\end{equation}
%where $c_{n+1} = 4n/(n-1)$, $p = (n+2)/(n-1)$.  
Let $v = u-1$ and let  
$$
\hS = \frac{d-2}{4(d-1)} \left ( S + d(d-1) \right ).
$$
Then the Yamabe equation takes the form 
\begin{equation}\label{eq:Yam-v}
- \Delta v + d v + \hS v = - \hS - \calF(v)
\end{equation} 
where 
$$
\calF(v) = \frac{d(d-2)}{4} \left [ (1+v)^{\frac{d+2}{d-2}} -1 -
  \frac{d+2}{d-2} v \right ] 
$$
In particular $\calF(v) = O(v^2)$. 
A straightforward application
of the maximum principle shows that since $S[g] \geq - d(d-1)$, we have  $v
\leq 0$ and hence 
$$
u \leq 1  \,.
$$ 
Linearizing the Yamabe equation around $u = 1$, we
obtain the equation 
$$
- \Delta \bu + d \bu + \hS \bu = 0   \,.
$$
The indicial exponents of this equation are  $-1, d$.
It follows that the solution to the Yamabe equation is of the form 
$u = 1 + v
$
with 
$
v = v_{d,1} t^d \log t + v_d t^d + \text{higher order}. 
$
However since by Equation \eqref{eq:scal},  $\hS = O(t^{d+1})$,
it follows \cite{and:chru}  that $v_{d,1} = 0$, and in
fact $v$ is smooth up to boundary, with 
$$
v = v_d t^d + \text{ higher order} .
$$ 

Let $L = - \Delta + d$. Equation (\ref{eq:Yam-v}) takes the form 
\begin{equation}\label{eq:Lv=f} 
L v =  f
\end{equation} 
with $f$ given by 
$$
f = - \hS u - \calF(u-1)  \,.
$$
In particular, $f \leq 0$ and $f \ne 0$ except when $\hS = 0$. 
Let $L_t$ be the operator defined by 
$$
L_t u = - \sinh^2(t) \partial_t^2 u + (d-2) \sinh(t)\cosh(t) \partial_t u 
+ d u  \,.
$$
We have 
$$
L u = L_t u  - \sinh^2(t) \partial_t \sqrt{\det h} \partial_t u -
\sinh^2(t) \Delta_h u  \,,
$$
where $\Delta_h$ is the Laplacian on $S^{d-1}$ with respect to the
metric $h(t,\cdot) = h_0 + t^d \gamma$. In particular $\Delta_h$ involves only
$y^A$-derivatives. 

We now introduce a function $w$ which will be used as a supersolution, in order
to control the leading order term in $v$. Let $w = - t^d(1+d t)$. 
\begin{lem} \label{lem:Lt-est} There exists constants 
$t_1 = t_1(d) > 0$, $A = A(d) > 0$  such that 
$$
L_t w > A t^{d+1}, \quad \text{ for } 0 < t < t_*
$$
\end{lem} 
\begin{proof}  Using $w = - t^d(1+d t)$, we obtain,
\begin{align*}
L_t(w)&=\sinh^2(t)[d(d-1)t^{d-2}+d^2(d+1)t^{d-1}]
\\&\qquad -(d-2)\sinh(t)\cosh(t)[dt^{d-1}+d(d+1)t^{d}] -d(t^d+dt^{d+1})
\\&=[d(d-1)t^d+d^2(d+1)t^{d+1}]-(d-2)[dt^{d}+d(d+1)t^{d+1}]
\\&\qquad - [dt^{d}+d^2t^{d+1}] + O(t^{d+2})\\
&=d(d+2)t^{d+1}+O(t^{d+2}) \,,
\end{align*}
and the lemma follows.
\end{proof}
We have, 
$$
L w = L_t w - \sinh^2(t) \partial_t \sqrt{\det h} \partial_t w  \,.
$$
For metrics of the form we are considering, $\tr_{h} \partial_t h =
O(t^{d-1})$. 
Hence by Lemma \ref{lem:Lt-est}, we have 
$$
L w >  A t^{d+1} - C t^{2d}
$$
where $C = C(d,\gamma)$. This means there exists $t_2 > 0$, $t_2 =
t_2(d,\gamma)$ such that 
$$
L w > 0, \quad \text{ for } 0 < t < t_2 \,.
$$

\begin{lem} \label{lem:green} 
Let $(M,g)$ be asymptotically hyperbolic, so that \eqref{rhometric} and
\eqref{rhometric2} hold.
Let $f$ be a function on $(M,g)$, and assume that $f$ is smooth up to
$\partial \tM$, with fall off $f = O(t^{d+1})$. Let $L = - \Delta + d$, and
let $v$ be the unique
solution to 
$$
L v = f
$$
with $v = O(t^d)$. Then $v = v_d t^d + t^{d+1} J$ with $v_d = v_d(y)$ smooth
on $\partial \tM$ and 
$J$ smooth up to
$\partial \tM$. If $f \leq 0$, $f \ne 0$, then  
$v_d < 0$. 
\end{lem} 
\begin{proof} 
Let $\bar v_a  =  a w$ and let $v$ be as in (\ref{eq:Lv=f}). We have $f \leq
0$, and hence by the strong maximum prinicple, $v < 0$ in the interior of
$\tilde M$. 
It follows that there is $\eps > 0$, $t_2 > t_3 > 0$, so that 
$$
\sup_{y \in S^{d-1}} v(t_3,y) < - \eps 
$$
For each $a \geq 0$, $\bar v_a$ is a supersolution to $L$ in the region $0 < t <
t_3$ and for $0 \leq a \leq a_*$, we have that $\bar v_a(t_3) > v(t_3,y)$ for
$y \in S^{d-1}$. Further, we clearly have $\bar v_a(0) = v(0,y) = 0$ for $y \in
S^{d-1}$. It follows from the maximum principle 
that for small $a$, $\bar v_a > v$ in the region $0 < t < t_3$. Fix an $a$
with this property. 
Since $\bar v_a = - a t^d + O(t^{d+1})$, dividing the inequality,
$ v \leq \bar v_a$, by $t^d$ and letting $t \searrow 0$ gives 
$v_d \le - a$. 
\end{proof}

We are now ready to state the following analogue of a well-known
result in the asymptotically flat setting (cf., \cite{niall}).

\begin{prop}\label{prop:reduce}
Let $(M,g)$ be asymptotically hyperbolic in the sense of Definition
\ref{def:ah}, with scalar curvature $S[g] \ge -d(d-1)$,  and with strict inequality
somewhere.  Then there exists a conformally related metric
$\hat g$ such that
\ben
\item $(M,\hat g)$ is asymptotically hyperbolic,
\item $S[\hat g] = -d(d-1)$, and
\item $\mu[\hat g] < \mu[g]$,
\een
where  $\mu[g]$,  $\mu[\hat g]$ are the mass aspect functions of $(M,g)$,
$(M,\hat g)$, respectively.
\end{prop}

\begin{proof} 
Lemma \ref{lem:green},  together with the discussion prior to Lemma \ref{lem:Lt-est},
shows that the solution $u$ to the Yamabe equation 
is of the form 
$$
u = 1 + u_d t^d + t^{d+1} J
$$
with $u_d  = u_d(y) < 0$, and with $J$ smooth up to $\partial \tM$. 
Then, after a change of coordinates, $\hme = u^{4/(d-2)} g$ can be brought into
the form,
$$
\hme  = \sinh^{-2}(t) (dt^2 + h_0 + t^d (\gamma+ \frac{4}{d-2}(1 + \frac1{d}) u_d
h_0) + t^{d+1} z )
$$
where $z= z_{ij}dx^idx^j$ is smooth up to $\partial \tM$. 

By Lemma \ref{lem:gauge}, after a change of conformal gauge, we have  
$$
\hme = \sinh^{-2} (\that) 
( d\that^2 + h_0 + \that^d (\hgamma+ \frac{4}{d-2}(1 + \frac1{d}) u_d h_0)  )
$$
where $\hgamma = \hgamma(\that,y)_{AB} dy^A dy^B$ 
is smooth up to $\partial \tM$ and $\hgamma(0,y) = \gamma(0,y)$. 
It follows from the above that the mass aspect functions
satisfy,
\beq
\mu[\hat g] = \mu[g] + \frac{4(d-1)}{d-2}(1 + \frac1{d}) u_d  < \mu[g] \,,
\eeq
since $u_d < 0$.
\end{proof}

For the purpose of establishing Theorem \ref{thm:rigid2}, we need the following
immediate consequence of Propositions  \ref{prop:aspect} and \ref{prop:reduce}. 

\begin{cor}\label{cor:confdef} 
Let $(M,g)$ be as in Theorem \ref{thm:rigid2}; in particular, assume
$\mu[g] = 0$. Then $g$ has constant scalar
curvature $S[g] = - d(d-1)$. 
\end{cor} 

\begin{proof}
Suppose $S[g] > - d(d-1)$ somewhere.  Then, by Proposition \ref{prop:reduce},
there exists a conformally related metric $\hat g$ such that $(M,\hat g)$ is 
asymptotically hyperbolic, $S[\hat g] = -d(d-1)$, and $\mu[\hat g] < \mu[g] = 0$.
But this directly contradicts Proposition \ref{prop:aspect}.  
\end{proof}

\subsubsection{Deforming the metric} 
Now we will show that if $g$ has constant scalar curvature $S = - d(d-1)$ and
vanishing mass aspect function, then it is Einstein, 
$\Ric_g = -(d-1) g$. Thus, let $(M,g)$ be as in Theorem \ref{thm:rigid2},
and 
assume 
$S[g] = -d(d-1)$. 

Let 
$$
\hRic = \Ric - \frac{S}{d} g
$$
denote the traceless part of $\Ric$. Note that since $g$ has constant scalar
curvature $\Ric$ and $\hRic$ have vanishing divergence.

For the subsequent analysis, we shall need detailed information about the asymptotic
behavior of $\hRic$. 
\begin{lem}\label{lem:Ricasympt} 
Let $(M,g)$ be as in Theorem \ref{thm:rigid2} (so that  \eqref{rhometric} and
\eqref{rhometric2} hold, and the mass aspect vanishes).
Then 
$$
\hRic = - \frac{d}{2} t^{d-2} \gamma + t^{d-1} z
$$
where $z = z_{ij} dx^i dx^j$. 
\end{lem} 
\begin{proof} 
Let 
$\tRic, \tnabla$ denote the Ricci tensor and covariant derivative defined
with respect to $\tme$. 
The conformal transformation formula for
Ricci curvature is 
\begin{equation}\label{eq:confRic} 
\Ric_{ij} = \tRic_{ij} + \rho^{-1} [ (d-2) \tnabla_i \tnabla_j \rho +
  \tnabla^l \tnabla_l \rho \tme_{ij} ] - (d-1) \rho^{-2} \tnabla_l \rho
\tnabla^l \rho 
\tme_{ij} 
\end{equation} 
where in the right hand side, indices are raised with $\tme$. 

Note that if $\gamma = 0$, then $\tRic = (d-2) h_0$ and the only nonvanishing
terms in the formula for $\Ric$ are
\begin{align*} 
\Ric_{tt} &= (d-1) \sinh^{-2}(t) \\ 
\Ric_{AB} &= (d-1) \sinh^{-2}(t) h_{0\, AB} 
\end{align*} 
Now we consider the case with nonvanishing $\gamma$, but with vanishing mass
aspect function, i.e.
$\mu = h_0^{AB} \gamma_{AB}
\big{|}_{t=0} = 0$. 

Let $K_{ij} = \half \partial_t \tme_{ij}$. Then the 
nonvanishing terms in $K$ are $K_{AB} = \half d t^{d-1} \gamma_{AB} +
O(t^d)$. Let 
$\hpara, \nablapara, \Ricpara$ denote the induced
metric, covariant derivative and Ricci tensor on the level sets $M_t$ of
$t$. 
We use coordinates $y^A$ on these level sets, and raise and lower indices
with $\hpara$. Note that $\hpara = h_0 + t^d \gamma$, and hence, since
$\Ricpara$ involves no $t$-derivatives, we have by Taylor's
theorem, 
$$
\Ricpara = \Ric[h_0] + O(t^d) = (d-2) h_0 + O(t^d)
$$

We have from the Gauss, Codazzi, and second variation equations  
\begin{align*} 
\tRic_{tt} &= - \hpara^{AB} \partial_t K_{AB} + K_{AC} K^C{}_B \\
&= - \half d(d-1) t^{d-2} \hpara^{AB} \gamma_{AB} + O(t^{d-1}) \\ 
\intertext{which using $\mu = 0$ gives,}
&= O(t^{d-1}) \,, \\ 
\tRic_{tA} &= \nablapara^B K_{BA} - \nablapara_A (\hpara^{BC} K_{BC} ) \\ 
&= O(t^d)  \,,\\ 
\tRic_{AB} &=  \Ricpara_{AB} - \partial_t K_{AB}  + 2 K_{AC} K^C{}_B - K_{AB}
\hpara^{CD} K_{CD} \\ 
&= (d-2) h_0 - \half d(d-1) t^{d-2} \gamma_{AB}  + O(t^{d-1})   \,.
\end{align*} 
Thus we have 
$$
\tRic = (d-2) h_0 - \half d(d-1) t^{d-2} \gamma + O(t^{d-1})
$$
We next consider the remaining terms 
$$
B_{ij} = \rho^{-1} [ (d-2) \tnabla_i \tnabla_j \rho +
  \tnabla^l \tnabla_l \rho \tme_{ij} ] - (d-1) \rho^{-2} \tnabla_l \rho
\tnabla^l \rho 
\tme_{ij}  \,.
$$
Recall that since we are in a Gauss
  foliation, the only non-vanishing terms in $\tGamma^t_{ij}$ are 
$$
\tGamma^t_{AB} = - K_{AB} \,.
$$
We have 
\begin{align*} 
B_{tt} &= - (d-1)\sinh^{-2}(t) \,, \\
B_{tA} &= 0 \, ,\\
B_{AB} &= \sinh^{-1}(t) [ (d-2) K_{AB} \cosh(t)  
+ \sinh(t) \hpara_{AB} ]  \\
&\quad -
(d-1) \sinh^{-2}(t) \cosh^2(t) \hpara_{AB} \,, \\
&= - (d-1) \sinh^{-2}(t)  h_{0\, AB} + \half d(d-2)t^{d-2} \gamma_{AB} + O(t^{d-1})  \,.
\end{align*} 
This shows that 
$$
\Ric[g] = - (d-1) g - \frac{d}{2} t^{d-2} \gamma + O(t^{d-1}) 
$$
which gives the Lemma. 
\end{proof} 

Consider the curve, 
\begin{equation}\label{eq:lambdadef} 
\lambda_s = u_s^{4/(d-2)} g_s
\end{equation} 
where, for $s$ small, $g_s$ is the smooth curve of metrics, $g_s = g -s\hRic[g]$,
%$\partial_s g_s \big{|}_{s=0} = - \tRic[g]$
and $u_s$ is the conformal
factor such that $S[\lambda_s] = -d(d-1)$.  Note that $u_0 = 1$.  By Lemmas
\ref{lem:gauge}, \ref{lem:Ricasympt}, and our earlier discussion
on the asymptotic form of solutions to the Yamabe equation, 
$\lambda_s$ is asymptotically hyperbolic in the sense of Definition~\ref{def:ah}.
Let $\mu_s$ denote the
mass aspect function of $\lambda_s$.

Let $\bar u = \frac{\d u_s}{\d s}|_{s=0}$.  Then, by differentiating the Yamabe 
equation \eqref{eq:yam}, with  $u = u_s$ and $g = g_s$, with respect to the parameter $s$, we obtain the equation,
$$
-\Delta \bu + d \bu = - |\hRic|^2  \,.
$$
By Lemma \ref{lem:green}, we have 
$$
\bu =  \bu_d t^d + O(t^{d+1})  \,,
$$
with $\bu_d  = \bu_d(y) < 0$ if $\hRic \ne 0$.
\begin{lem} \label{lem:dmu}  
$$
\partial_s \mu_s\big{|}_{s=0} = \frac{4(d-1)}{d-2}(1 + \frac1{d}) \bu_d  \,.
$$
\end{lem} 
\begin{proof} 
Clearly, $\alpha = \partial_s \mu_s
\big{|}_{s=0}$ is of the form $\alpha = \alpha_{\hRic} + \alpha_u$, where 
$$
\alpha_{\hRic} = \partial_s \mu (g_s) \big{|}_{s=0} 
$$
and 
$$
\alpha_u = \partial_s \mu(u_s^{4/(d-2)} g) \big{|}_{s=0}  \,.
$$
Let $g_s = \sinh^{-2}(t) \tme_s$. In order to determine the $s$-dependence of
the mass aspect function we consider $\tme_s$.  
It follows from Lemma \ref{lem:Ricasympt} that 
$$
\sinh^2(t)\partial_s \tme_s =  \frac{d}{2} t^{d}\gamma +
O(t^{d+1})   \,.
$$
Since by assumption $\tr_{h_0} \gamma \big{|}_{\partial M} = 0$, we have 
$\alpha_{\hRic} = 0$. It follows that the first order change in the mass
aspect function of $\lambda_s$ is given by $\alpha_u$, 
which clearly is determined by the first
order change in the conformal factor $u_s$. The result follows. 
\end{proof} 

We are now ready to complete the proof of Theorem \ref{thm:rigid2}. 
\begin{proof}[Proof of Theorem \ref{thm:rigid2}] 
Let $(M,g)$ be as in Theorem \ref{thm:rigid2}. Recall  $d =
n+1$. By Corollary
\ref{cor:confdef}, $S[g] = -d(d-1)$. Suppose that $g$ is not Einstein,
i.e. $\hRic \ne 0$. 
Let $\lambda_s = u_s^{4/(d-2)} g_s$
as in (\ref{eq:lambdadef}).  As previously observed, $\l_s$ is asymptotically
hyperbolic with  scalar curvature $S[\lambda_s] = -d(d-1)$, and with mass aspect $\mu_s$.  
By Lemma \ref{lem:dmu},  $\partial_s \mu_s|_{s=0} < 0$, and hence for small $s >
0$, $\lambda_s$ has negative mass aspect function.  But in view of 
Proposition \ref{prop:aspect}, $\mu_s < 0$ gives a contradiction, and hence it must hold that $\hRic = 0$. 
We can now apply the rigidity result of Qing \cite{qing} (see also
\cite{bonini:miao:qing:ricci,cai:qing}) 
to conclude that in fact $(M,g)$
is isometric to hyperbolic space. This concludes the proof of the positive
mass theorem in the case of vanishing mass aspect function. 
\end{proof} 

Naturally, it would be desirable to find a way to remove the sign condition
on the mass aspect from our positive mass result.  
Within the context of the approach taken in this paper, one possible way
to accomplish this would be to extend the results of Corvino-Schoen \cite{corvino,
corvino:schoen}
and Chru\'sciel-Delay \cite{chru:delay} on initial data
deformations to the asymptotically hyperbolic setting.  Starting from Proposition
\ref{prop:reduce}, the aim would be to deform the time-symmetric initial data to be 
exactly Schwarzschild-AdS outside a compact set, without changing the scalar 
curvature and the sign of the mass.   Starting from Schwarzschild-AdS with
mass $m < 0$,  the deformation result of section \ref{deform} is 
then easily proved.

\section*{Acknowledgements}
We thank Piotr Chru\'sciel for comments and discussion.
This work was initiated at the Conference on Mathematical Aspects of 
Gravitation at 
the Mathematical Research Institute in Oberwolfach in 2003.  We thank the
Institute for its hospitality and support. 
This work was supported in part by NSF grants DMS-0407732 and
DMS-0405906.

\bibliographystyle{amsplain}
\bibliography{pmass}

\providecommand{\bysame}{\leavevmode\hbox to3em{\hrulefill}\thinspace}
\providecommand{\MR}{\relax\ifhmode\unskip\space\fi MR }
\providecommand{\MRhref}[2]{%
  \href{http://www.ams.org/mathscinet-getitem?mr=#1}{#2}
}
\providecommand{\href}[2]{#2}
\begin{thebibliography}{10}

\bibitem{and:chru}
Lars Andersson and Piotr~T. Chru{\'s}ciel, \emph{Solutions of the constraint
  equations in general relativity satisfying ``hyperboloidal boundary
  conditions''}, Dissertationes Math. (Rozprawy Mat.) \textbf{355} (1996), 100.
  \MR{MR1405962 (97e:58217)}

\bibitem{AD}
Lars Andersson and Mattias Dahl, \emph{Scalar curvature rigidity for
  asymptotically locally hyperbolic manifolds}, Ann. Global Anal. Geom.
  \textbf{16} (1998), no.~1, 1--27.

\bibitem{aviles:mcowen}
Patricio Aviles and Robert McOwen, \emph{Conformal deformations of complete
  manifolds with negative curvature}, J. Differential Geom. \textbf{21} (1985),
  no.~2, 269--281. \MR{MR816672 (87e:53058)}

\bibitem{bonini:miao:qing:ricci}
Vincent Bonini, Pengzi Miao, and Jie Qing, \emph{Ricci curvature rigidity for
  weakly asymptotically hyperbolic manifolds}, Comm. Anal. Geom. \textbf{14}
  (2006), no.~3, 603--612.

\bibitem{cai}
Mingliang Cai, \emph{Volume minimizing hypersurfaces in manifolds of
  nonnegative scalar curvature}, Minimal surfaces, geometric analysis and
  symplectic geometry (Baltimore, MD, 1999), Adv. Stud. Pure Math., vol.~34,
  Math. Soc. Japan, Tokyo, 2002, pp.~1--7.

\bibitem{CG}
Mingliang Cai and Gregory~J. Galloway, \emph{Rigidity of area minimizing tori
  in 3-manifolds of nonnegative scalar curvature}, Comm. Anal. Geom. \textbf{8}
  (2000), no.~3, 565--573.

\bibitem{cai:qing}
Mingliang Cai and Jie Qing, \emph{On the rigidity of {A}d{S} spacetime},
  preprint (2006).

\bibitem{christ:loh}
U.~Christ and J.~Lohkamp, \emph{Singular minimal hypersurfaces and scalar
  curvature}, math.DG/0609338 (2006).

\bibitem{chru:delay}
Piotr~T. Chru{\'s}ciel and Erwann Delay, \emph{On mapping properties of the
  general relativistic constraints operator in weighted function spaces, with
  applications}, M\'em. Soc. Math. Fr. (N.S.) (2003), no.~94, vi+103.
  \MR{MR2031583 (2005f:83008)}

\bibitem{ChH}
Piotr~T. Chru{\'s}ciel and Marc Herzlich, \emph{The mass of asymptotically
  hyperbolic {R}iemannian manifolds}, Pacific J. Math. \textbf{212} (2003),
  no.~2, 231--264.

\bibitem{corvino}
Justin Corvino, \emph{Scalar curvature deformation and a gluing construction
  for the {E}instein constraint equations}, Comm. Math. Phys. \textbf{214}
  (2000), no.~1, 137--189. \MR{MR1794269 (2002b:53050)}

\bibitem{corvino:schoen}
Justin Corvino and Richard~M. Schoen, \emph{On the asymptotics for the vacuum
  {E}instein constraint equations}, J. Differential Geom. \textbf{73} (2006),
  no.~2, 185--217. \MR{MR2225517}

\bibitem{delay}
Erwann Delay, \emph{Analyse pr\'ecis\'ee d'\'equations semi-lin\'eaires
  elliptiques sur l'espace hyperbolique et application \`a la courbure scalaire
  conforme}, Bull. Soc. Math. France \textbf{125} (1997), no.~3, 345--381.

\bibitem{federer}
Herbert Federer, \emph{Geometric measure theory}, Die Grundlehren der
  mathematischen Wissenschaften, Band 153, Springer-Verlag New York Inc., New
  York, 1969. \MR{MR0257325 (41 \#1976)}

\bibitem{gibbonsetal}
G.~W. Gibbons, S.~W. Hawking, Gary~T. Horowitz, and Malcolm~J. Perry,
  \emph{Positive mass theorems for black holes}, Comm. Math. Phys. \textbf{88}
  (1983), no.~3, 295--308. \MR{MR701918 (84k:83015)}

\bibitem{KW}
Jerry~L. Kazdan and F.~W. Warner, \emph{Existence and conformal deformation of
  metrics with prescribed {G}aussian and scalar curvatures}, Ann. of Math. (2)
  \textbf{101} (1975), 317--331.

\bibitem{loh2}
J.~Lohkamp, \emph{The higher dimensional positive mass theorem {I}},
  math.DG/0608795 (2006).

\bibitem{loh}
Joachim Lohkamp, \emph{Scalar curvature and hammocks}, Math. Ann. \textbf{313}
  (1999), no.~3, 385--407.

\bibitem{minoo}
Maung Min-Oo, \emph{Scalar curvature rigidity of asymptotically hyperbolic spin
  manifolds}, Math. Ann. \textbf{285} (1989), no.~4, 527--539.

\bibitem{niall}
Niall \'O~Murchadha and James~W. York~Jr., \emph{Gravitational energy}, Phys.
  Rev. D \textbf{10} (1974), no.~8, 2345--2357.

\bibitem{qing}
Jie Qing, \emph{On the rigidity for conformally compact {E}instein manifolds},
  Int. Math. Res. Not. (2003), no.~21, 1141--1153.

\bibitem{SY}
R.~Schoen and S.~T. Yau, \emph{On the structure of manifolds with positive
  scalar curvature}, Manuscripta Math. \textbf{28} (1979), no.~1-3, 159--183.

\bibitem{SY1}
Richard Schoen and Shing~Tung Yau, \emph{On the proof of the positive mass
  conjecture in general relativity}, Comm. Math. Phys. \textbf{65} (1979),
  no.~1, 45--76.

\bibitem{shi}
Yuguang Shi and Gang Tian, \emph{Rigidity of asymptotically hyperbolic
  manifolds}, Comm. Math. Phys. \textbf{259} (2005), no.~3, 545--559.

\bibitem{spivak}
Michael Spivak, \emph{A comprehensive introduction to differential geometry},
  second ed., Publish or Perish Inc., Wilmington, Del., 1979.

\bibitem{wang}
Xiaodong Wang, \emph{The mass of asymptotically hyperbolic manifolds}, J.
  Differential Geom. \textbf{57} (2001), no.~2, 273--299.

\bibitem{WY}
Edward Witten and S.-T. Yau, \emph{Connectedness of the boundary in the
  {A}d{S}/{CFT} correspondence}, Adv. Theor. Math. Phys. \textbf{3} (1999),
  no.~6, 1635--1655 (2000).

\bibitem{yau}
Shing~Tung Yau, \emph{Geometry of three manifolds and existence of black hole
  due to boundary effect}, Adv. Theor. Math. Phys. \textbf{5} (2001), no.~4,
  755--767.

\end{thebibliography}

\end{document}